\documentclass[11pt]{article}
\usepackage{amsfonts}
\usepackage{eucal}
\usepackage{amssymb,epsfig,latexsym,color}
\usepackage{graphicx,amsmath,amssymb,latexsym,psfrag}
\usepackage{graphics,graphicx,comment}
\usepackage{multirow}
\newcommand{\R}{\mathbb{R}}

\newcommand{\EE}{\mathbb{E}}

\newcommand\Top{\rule{0pt}{3.6ex}}       
\newcommand\Bottom{\rule[-2.2ex]{0pt}{0pt}} 
\newtheorem{Theorem}{Theorem}[section]
\newtheorem{Proposition}[Theorem]{Proposition}

\newtheorem{Lemma}[Theorem]{Lemma}
\newtheorem{Remark}[Theorem]{Remark}

\def\ind{{\rm 1\hspace{-0.90ex}1}}

\begin{document}
\title{A large deviation approach to super-critical bootstrap percolation on the random graph $G_{n,p}$}
\author{Giovanni
Luca Torrisi\thanks{Istituto per le Applicazioni del Calcolo, CNR, Roma, Italy.
e-mail: \tt{giovanniluca.torrisi@cnr.it}} \and Michele Garetto\thanks{Dipartimento di Informatica,
Universit\`{a} di Torino, Italy. e-mail: \tt{michele.garetto@unito.it}}\and Emilio Leonardi\thanks{Dipartimento di Elettronica,
Politecnico di Torino, Italy. e-mail: \tt{emilio.leonardi@polito.it}}
\thanks{{\bf Complete technical report version of a paper submitted to {\em Stochastic Processes and Their Applications}}
}}
\date{}
\maketitle
\begin{abstract}

We consider the Erd\"{o}s--R\'{e}nyi random graph $G_{n,p}$ and we analyze the simple
irreversible epidemic process on the graph, known in the literature as bootstrap percolation.
We give a quantitative version of some results
by Janson et al. (2012), providing a fine asymptotic analysis
of the final size $A_n^*$ of active nodes, under a suitable
super-critical regime. More specifically, we establish large deviation principles for the sequence of random variables
$\{\frac{n- A_n^*}{f(n)}\}_{n\geq 1}$ with explicit rate functions and allowing the scaling function $f$ to vary in the
widest possible range.

\end{abstract}

\noindent\emph{MSC 2010 Subject Classification}: 05C80, 60K35, 60F10.\\
\noindent\emph{Keywords}: Bootstrap Percolation, Large Deviations, Random Graphs.\\
\section{Introduction}

Bootstrap percolation on a graph is a simple activation process
that starts with a given number of initially active nodes (seeds) and evolves as follows.
An inactive node that has at least $r\geq 2$ active neighbors becomes active, and remains so forever.
The process stops when no more nodes become active.

Bootstrap percolation has a rich history and was initially investigated
on regular structures, see e.g. \cite{adler} for a survey. The study of bootstrap percolation
on lattices and grids can be explained by its origin
in the area of statistical physics.
Bootstrap percolation on a lattice was introduced
in \cite{chalupa} and further studied in \cite{schon}. Deep results for the bootstrap percolation process
over finite grids (in two dimensions or more)
were obtained by several authors \cite{AL,BP,BBM,BBDM,CC,CM,GHM,H}. Bootstrap percolation
on the hypercube and trees was investigated in \cite{bollobas} and \cite{BPP}, respectively.
We cite \cite{FSS} and \cite{M} for applications of bootstrap percolation to the Ising model.

More recently, the bootstrap percolation process has been investigated in the context of
random graphs. This is partly motivated by the increasing interest
in dynamical processes taking place over large-scale complex systems
such as technological, biological and social networks
whose irregular structure is better captured by random graphs models (see
\cite{masso} for a comprehensive introduction to epidemics in complex networks).
For example, in the case of social networks, bootstrap percolation may serve as a primitive
model for the spread of ideas, rumors and trends among individuals. Indeed, in this context
one can assume that a person will adopt an idea after receiving sufficient influence by
friends who have already adopted it \cite{kempe,munik,watts}.

Bootstrap percolation on random regular graphs was studied in \cite{pittel} and \cite{J}. This analysis
was extended in \cite{amini1} to random graphs with given vertex degrees (configuration model). A bootstrap percolation model where edges are activated rather than nodes
was introduced in \cite{bollobas1} and recently analyzed
in \cite{balogh-edges}. In \cite{geometric} the authors derived critical thresholds for the bootstrap percolation
process in random geometric graphs. Bootstrap percolation on random graphs
was investigated also from an algorithmic perspective
\cite{soda,kempe}, with the goal of identifying the set of seeds that maximizes the final size.

In the seminal paper \cite{JLTV}, Janson, Luczak, Turova and Vallier provided a detailed analysis
of the bootstrap percolation process on the Erd\"{o}s-R\'{e}nyi random graph $G_{n,p_n}$,
i.e., the random graph on the set of nodes $\{1,\ldots,n\}$ where any two nodes are connected with
probability $p_n$, independently of all other node pairs. In \cite{JLTV} the authors assume that $a_n$ seeds are initially chosen uniformly at random
among the nodes and, under suitable assumptions which imply a sub-linear growth of the number of seeds,
proved the existence of a sharp phase transition.
Roughly speaking, below a critical number of seeds $a_c^{(n)}$, whose value
is available in closed form, the process essentially does not evolve, reaching, as $n\to\infty$, a
final size of active nodes $A_n^*$ which is of the same order as $a_n$
(sub-critical case), i.e., in mathematical terms, $A_n^*/a_n$ converges in probability to a suitable positive constant.
Instead, above the critical number $a_c^{(n)}$, the process percolates
through the entire random graph, reaching, as $n\to\infty$,
a final size of active nodes which is of the same order as $n$
(super-critical case), i.e., in mathematical terms, $A_n^*/n$ converges to $1$ in probability.

In \cite{kang} the results of \cite{JLTV} were extended to $k$-uniform
random hypergraphs.
Bootstrap percolation on random graphs obtained by combining $G_{n,p_n}$ with a regular lattice
was investigated in \cite{turova}. We mention also the recent work
\cite{holmgren}, where the authors studied the so-called majority bootstrap percolation on $G_{n,p_n}$,
according to which nodes become active when the number of their active neighbors exceeds the number of inactive neighbors.

In this paper we consider the super-critical regime of bootstrap percolation on $G_{n,p_n}$ and provide a
deeper investigation of the results in \cite{JLTV}; more specifically, we prove
large deviation estimates for $A_n^*$.
Roughly speaking, for any Borel set $B$ and various scaling functions $f$ and speed functions $v$,
we establish asymptotic estimates of the form
\begin{equation}\label{eq:ldp}
P\left(\frac{n-A_n^*}{f(n)}\in B\right)\approx\mathrm{e}^{-\mathcal{I}(B)v(n)},\quad\text{as $n\to\infty$}
\end{equation}
where the quantity $\mathcal{I}(B)$ is explicitly given and the approximation is in the sense of large deviations
(see Section \ref{par:ldp} and Theorems \ref{ldp1}, \ref{ldp2}, \ref{ldp5}, \ref{ldp3}, \ref{ldp4}).
An estimate as \eqref{eq:ldp} provides a natural quantitative version of the corresponding weak law of large numbers in \cite{JLTV},
in the following sense. On the one hand, the weak law of large numbers determines the most probable value of $(n-A_n^*)/f(n)$, as $n\to\infty$;
on the other hand, the corresponding large deviation principle provides accurate estimates for the probability that $(n-A_n^*)/f(n)$
deviates from its most probable value. We postpone to Section \ref{sec:disc} an informal discussion of our results, and now proceed to describing the strategy of our proofs.

All the large deviation principles obtained in this paper are proved in the following way. Firstly,
we characterize the asymptotic behavior of $\log P((n-A_n^*)/f(n)>\varepsilon)$, with $\varepsilon>0$, as $n$ grows large (see Propositions \ref{ldp:supcrit}, \ref{ldp:supcritbis},
\ref{ldp:supcritpenta}, \ref{ldp:supcritquater}, \ref{ldp:supcritter}).
Secondly, we provide large deviation principles by combining such tail asymptotics with elementary topological considerations which allow us to study the asymptotic behavior of
$\log P((n-A_n^*)/f(n)\in B)$, as $n$ grows large, for any Borel set $B$ (see the proofs of Theorems \ref{ldp1}, \ref{ldp2}, \ref{ldp5}, \ref{ldp3}, \ref{ldp4}).

From a technical point of view, the first step is certainly the core of this paper. Its proof is based on a fine analysis of the set
$\{(n-A_n^*)/f(n)>\varepsilon\}$, as $n\to\infty$, aimed at determining the dominant event which characterizes the asymptotic behavior of the tail probability
on a logarithmic scale.
Basically, the set $\{(n-A_n^*)/f(n)>\varepsilon\}$ is re-written as the union of suitable events whose probabilities are estimated by
exploiting the \lq\lq binomial structure" of the bootstrap percolation process discovered in \cite{JLTV} (see also Subsection \ref{Sec:bin}).
Such estimates are obtained via concentration inequalities for the binomial distribution and other tools from the theory of large deviations.
We remark that the techniques of this paper differ substantially from those adopted in the seminal paper \cite{JLTV}, where the
authors employ Doob's inequality to address a less general problem.

As a by-product of large deviations, we strengthen the results obtained in \cite{JLTV}
for the super-critical regime, providing strong laws of large numbers for the final size
of active nodes (see Theorems \ref{lln2} and \ref{lln5}).

Finally, we emphasize that our results can be used as a building block to analyze
the bootstrap percolation process on random graphs more general than $G_{n,p_n}$.
Indeed, the simple structure of an Erd\"os-R\'enyi random graph
can often be recognized in suitable sub-graphs of more complex networks. For example, in~\cite{amini2}
bootstrap percolation on random graphs with power-law degree has been studied by applying the
results in \cite{JLTV} to a properly defined sub-graph with a sufficiently large number of nodes of high degree.
Random graph models capturing the community structure observed in many realistic systems
(like stochastic block models \cite{newman}) can be potentially analyzed exploiting similar ideas,
i.e., by jointly applying available results for $G_{n,p_n}$ to proper sub-graphs having an
Erd\"os-R\'enyi structure. Tight exponential estimates on the convergence rate of
the bootstrap percolation process on $G_{n,p_n}$, such as those derived here by large deviation principles,
may be needed when the number of sub-graphs is unbounded.
At last we wish to remark that complementary results to ours  have been obtained in 
 \cite{subcritico},  where  large deviation bounds  for $A_n^*$ 
in the sub-critical regime have been derived.

The paper is organized as follows. In Section \ref{Sect:model} we collect some preliminaries. Specifically, we provide the formal definition
of the bootstrap percolation process on $G_{n,p_n}$, we introduce some notation and we recall the notion of large deviation principle.
The main results of the paper are stated in Section \ref{sec:results}, discussed in Section \ref{sec:disc}
and proved in Section \ref{sec:proofs}. We include an Appendix
which contains the derivations of some auxiliary asymptotic relations.

\section{Preliminaries}\label{Sect:model}

\subsection{The bootstrap percolation process on the random graph $G_{n,p_n}$}\label{Sec:bin}

We consider the bootstrap percolation process on $G_{n,p_n}$ starting with
an initial set $\mathcal{A}_n(0)\subseteq\{1,\ldots,n\}$
of active nodes (seeds) of cardinality $a_n$, which are chosen uniformly at random
among the nodes of the random graph. Nodes not belonging to $\mathcal{A}_n(0)$ are initially inactive.
An inactive node becomes active as soon as at least $r$
of its neighbors are active, where $r \geq 2$ is a given integer.
Seeds are declared to be active irrespective of the
state of their neighbors. Active nodes never become inactive and so the set of active nodes
grows monotonically.

The bootstrap percolation process naturally evolves through generations
of nodes which are sequentially activated. The first generation is composed by all nodes having at least $r$ seeds as neighbors. The second generation is formed by all nodes having at least $r$ neighbors among the seeds and the nodes belonging to the first generation, and so on.
The bootstrap percolation process stops when either all of the nodes are active or
an empty generation is obtained.

To analyze the final number of active nodes it is convenient to adopt a problem reformulation,
originally proposed in \cite{scalia}, according to which a single node is activated at a time
(note that, by so doing, we forget about the generations). Specifically,
we introduce a virtual discrete time $t\in\mathbb N:=\{1,2,\ldots\}$
and we assign a mark counter $M_{i}^{(n)}(t)$, $M_i^{(n)}(0):=0$, to each inactive node $i$.
At time $t=1$ we arbitrarily choose $u_1\in\mathcal{A}_n(0)$ and add one mark to all its neighbors.
We say that node $u_1$ has been \lq used', we define
$\mathcal{U}_n(1):=\{u_1\}$ and we update the set of active nodes setting $\mathcal{A}_n(1):=\mathcal{A}_n(0)\cup\Delta\mathcal{A}_n(1)$,
where $\Delta\mathcal{A}_n(1)$ is the set of inactive nodes that become active at time $t=1$.
We continue recursively: at a generic time $t\in\mathbb N$, we choose a
node $u_t\in\mathcal{A}_n(t-1)\setminus\mathcal{U}_n(t-1)$, i.e., an active node that has not yet been
used, we add a mark to all its neighbors, we define
$\mathcal{U}_n(t):=\mathcal{U}_n(t-1)\cup\{u_t\}=\{u_s\}_{1\leq s\leq t}$ and we
update the set of active nodes setting $\mathcal{A}_n(t):=\mathcal{A}_n(t-1)\cup\Delta\mathcal{A}_n(t)$,
where $\Delta\mathcal{A}_n(t)$ is the set of inactive nodes that become active at time $t$.
Note that $\Delta\mathcal{A}_n(t)=\emptyset$ if $t<r$. The bootstrap percolation process terminates when
$\mathcal{A}_n(t)=\mathcal{U}_n(t)$.

Let $T_n$ be the time at which the process stops, i.e.,
\[
T_n:=\min\{t\in\mathbb N:\,\,\mathcal{A}_n(t)=\mathcal{U}_n(t)\}.
\]
Let $A_n(t)$ and $U_n(t)$ be the cardinality of $\mathcal{A}_n(t)$ and $\mathcal{U}_n(t)$, respectively.
Since $\mathcal{U}_n(t)\subseteq\mathcal{A}_n(t)$
and $U_n(t)=t$, we have
\begin{equation}
\label{eq:T}
T_n=\min\{t\in\mathbb N:\,\,A_n(t)=t\}=\min\{t\in\mathbb N:\,\,A_n(t)\leq t\}\quad\text{and}\quad A_n^*=T_n,
\end{equation}
where $A_n^*:=A_n(T_n)$ is the final size of the set of active nodes. For later purposes, note that
\begin{equation}\label{eq:ratiopos}
a_n\leq A_n^*\leq n.
\end{equation}

We now introduce an alternative description of the random variable $A_n(t)$.
For $s\leq T_n$ and $i\notin\mathcal{U}_n(s)$, let $I_i^{(n)}(s)$ be the indicator
that there is an edge between node $u_s$ and node $i$, i.e., the indicator that $i$
gets a mark by the node used at time $s$.
It follows that the number of marks that $i\notin\mathcal{U}_n(t)$
has accumulated at time $t\le T_n$ is
\begin{equation}\label{eq:mit}
M_i^{(n)}(t)=\sum_{s=1}^{t}I_i^{(n)}(s).
\end{equation}
The random variables $I_i^{(n)}(s)$, $s\leq T_n$, $i\notin\mathcal{U}_n(t)$, are independent and Bernoulli distributed with mean $p_n$, hence
$M_i^{(n)}(t)$ has the same law of $\mathrm{Bin}(t,p_n)$, where $\mathrm{Bin}(m,p)$ denotes a random variable distributed according
to the binomial law with parameters $(m,p)$.  Furthermore, note that, for any $i\notin\mathcal{U}_n(t)$ and
$1\leq t\leq T_n$, we have $i\in\mathcal{A}_n(t)$ if and only if $M_i^{(n)}(t)\ge r$. We have defined the random marks $I_i^{(n)}(s)$
for $s\leq T_n$ and $i\notin\mathcal{U}_n(t)$, but, as noticed in \cite{JLTV},
it is possible to introduce additional, redundant 
random marks, which are independent and Bernoulli distributed with mean $p_n$, in such a way
that $I_i^{(n)}(s)$ is defined for all $i\in\{1,\ldots,n\}$ and $s\in\mathbb N$.
Such additional marks are added, for any $s\le T_n$,
to already active nodes and so they have no effect on the
underlying bootstrap percolation process.
The gain of this construction is that, for any $t\in\mathbb N$, we can consider
a sequence $\{M_i^{(n)}(t)\}_{1\leq i\leq n}$ of independent and identically distributed random variables
expressed by \eqref{eq:mit}.

We define
\[
Y^{(n)}_i:=\min\{t\in\mathbb N:\,\,M_i^{(n)}(t)\geq r\}
\]
and observe that  for every $i \notin  \mathcal A_n(0)$, if $Y_i^{(n)}\leq T_n$, then $Y_i^{(n)}$ is the time at which node $i$ becomes active.
We clearly have
\[
\mathcal{A}_n(t)=\mathcal{A}_n(0)\cup\{i\notin\mathcal{A}_n(0):\,\,M_i^{(n)}(t)\geq r\}, \quad  {t\le T_n},
\]
and so, defining
\[
S_n(t):=\sum_{i\notin\mathcal{A}_n(0)}\ind\{Y_i^{(n)}\leq t\}\quad\text{and}\quad A_n(0):=a_n\in\mathbb N,
\]
we get
\begin{equation}\label{eq:A(t)}
A_n(t)=a_n+S_n(t),\quad  {t\le T_n}.
\end{equation}
We note that $S_n(t)$ is distributed as $\mathrm{Bin}(n-a_n,\pi_n(t))$, where
\begin{equation}\label{eq:pi}
\pi_n(t):=P(Y^{(n)}_1\leq t)=P(\mathrm{Bin}(t,p_n)\geq r),\quad t\in\mathbb N\cup\{0\}.
\end{equation}
We remark that  in the following, with  a small abuse of notation,  we will estend the definitions 
of  $\pi_n(t)$ and  $S_n(t)$ also to $t>T_n$, as follows:  $\pi_n(t):= P(\mathrm{Bin}(t,p_n)>r)$ and  
 $S_n(t):=\mathrm{Bin}(n-a_n,\pi_n(t))$. 
Furthermore,  for $t\in\mathbb N$,
\begin{equation}\label{eq:dis2}
\{S_n(t)+a_n\leq t\}\subseteq\{S_n'(t)\leq t\},
\end{equation}
where
\begin{equation}\label{eq:snprime}
S_n'(t):=S_n(t)+\mathrm{Bin}(a_n,\pi_n(t))
\end{equation}
with $\mathrm{Bin}(a_n,\pi_n(t))$ independent of $S_n(t)$. We extend the definition of $S_n$, $\pi_n$ and $S_n'$ to $\mathbb R_+:=(0,\infty)$ by setting
$S_n(t):=S_n(\lfloor t\rfloor)$, $\pi_n(t):=\pi_n(\lfloor t\rfloor)$ and $S_n'(t):=S_n'(\lfloor t\rfloor)$, $t\in\mathbb R_+$.
\begin{equation}\label{eq:snprime}
\end{equation}
We extend the definition of $S_n$, $\pi_n$ and $S_n'$ to $\mathbb R_+:=(0,\infty)$ by setting
$S_n(t):=S_n(\lfloor t\rfloor)$, $\pi_n(t):=\pi_n(\lfloor t\rfloor)$ and $S_n'(t):=S_n'(\lfloor t\rfloor)$, $t\in\mathbb R_+$.
Hereafter, for $x\in\mathbb R$, we put
\[
\lfloor x\rfloor:=\max\{m\in\mathbb Z:\,\,m\leq x\}\quad\text{and}\quad\lceil x\rceil:=\min\{m\in\mathbb Z:\,\,m\geq x\}.
\]

\subsection{Further notation and assumptions}\label{sec:prel}

We shall use the following asymptotic notation. Let $f,g:\mathbb N\to\R$ be two functions. We write: $f(n)=o(g(n))$, or, equivalently, $f(n) \ll g(n)$,
if $\lim_{n\to\infty}\frac{f(n)}{g(n)}=0$; $f(n)=O(g(n))$ if there exist $c>0$, $n_0\in\mathbb N$: $|f(n)|\leq c|g(n)|$, for any $n\geq n_0$;
$f(n)\sim g(n)$ if there exists $c\in\R\setminus\{0\}$: $\lim_{n\to\infty}\frac{f(n)}{g(n)}=c$ and $f(n)\sim_e g(n)$ if $c=1$;
$f(n)\lesssim g(n)$ if either $f(n)\ll g(n)$ or $f(n)\sim g(n)$.
Unless otherwise specified, in this paper all the limits are taken as $n\to\infty$. We denote by $a\wedge b$ and $a\vee b$ the minimum and the maximum between $a,b\in\mathbb R$,
respectively.

As in \cite{JLTV}, throughout this paper we shall assume
\begin{equation}\label{eq:hyp}
1/(np_n)\to 0\quad\text{and}\quad p_n=o(n^{-1/r}),
\end{equation}
and we shall consider the following critical quantities, which allow us
to discriminate among different regimes:
\begin{eqnarray}
t_c^{(n)} & := & \left( \frac{(r-1)!}{n p_n^r} \right)^{1/(r-1)} \label{eq:tcjanson} \\
a_c^{(n)}  & := & \left(1 - \frac{1}{r}\right)t_c^{(n)} \label{eq:acjanson} \\
b_c^{(n)} & := &n\frac{(n p_n)^{r-1}}{(r-1)!}\mathrm{e}^{-n p_n}.  \label{eq:bcjanson}
\end{eqnarray}
As already mentioned in the Introduction, $a_c^{(n)}$ represents the critical number of seeds
associated to the phase transition between the sub-critical and the super-critical case.
In this paper we shall only consider the super-critical bootstrap percolation, i.e., we shall assume
\begin{equation}\label{hyp:supcrit}
a_n/a_c^{(n)}\to\alpha>1.
\end{equation}
The quantity $t_c^{(n)}$ represents the critical time associated
to $a_c^{(n)}$, while
$b_c^{(n)}$ has the same asymptotic behavior as the mean number of nodes with degree strictly less than $r$ and, as we shall see at the end of Subsection \ref{subsec:ldp},
characterizes different regimes for the final size of active nodes.

It is of rather immediate verification (see \cite{JLTV}) that, under \eqref{eq:hyp},
\begin{equation}\label{eq:ptcto0}
a_c^{(n)}\to+\infty,\quad a_c^{(n)}/n\to 0,\quad p_n a_c^{(n)}\to 0\quad\text{and}\quad b_c^{(n)}=o\left(\frac{a_c^{(n)}}{np_n }\right).
\end{equation}

\subsection{Large deviation principles}\label{par:ldp}

We say that a family of probability measures $\{\mu
_n\}_{n\in\mathbb N}$ on a topological space $(M,\mathcal{T}_M)$ obeys a large deviation principle (LDP) on $M$ with rate function
$I$ and speed $v$ if $I:M\rightarrow [0,\infty]$ is a lower semi-continuous function,
$v:\mathbb N\rightarrow (0,\infty)$ is a measurable function
which diverges to infinity, and the following
inequalities hold:
\[
\liminf _{n\rightarrow\infty}\frac{1}{v(n)}\log\mu_n(O)\geq-\inf_{x\in O}I(x),\quad\text{for every open set $O\subseteq M$}
\]
and
\[
\limsup_{n\rightarrow\infty}\frac{1}{v(n)}\log\mu_{n}(C)\leq-\inf
_{x\in C}I(x),\quad\text{for every closed set $C\subseteq M$.}
\]
Similarly, we say that a family of $M$-valued random variables
$\{V_n\}_{n\in\mathbb N}$ obeys an LDP on $M$ with rate function $I$ and speed $v$
if $\{\mu_n\}_{n\in\mathbb N}$, $\mu_n(\cdot):=P(V_n\in\cdot)$, obeys an LDP on $M$ with rate function $I$ and speed $v$.
We refer the reader to \cite{DZ} for an introduction to the theory of large deviations.

\section{Main results}\label{sec:results}

In this section we state our main results, referring the reader   to Section~\ref{sec:disc}  for
an informal discussion.

We define the following functions:
\begin{equation*}
h(x):=r^{-1}(\alpha(1-r^{-1})+x)^r\quad\text{and}\quad J(x):=\frac{r}{r-1}h(x)H\left(\frac{x}{h(x)}\right),\quad x\geq 0
\end{equation*}
where
\begin{equation}\label{eq:Hx}
H(x):=1-x+x\log x,\quad x\in\R_+,\quad H(0):=1,\quad H(x):=+\infty,\quad x\in\R_-:=(-\infty,0).
\end{equation}
We denote by $x_0$ the unique point of minimum over $[0,\infty)$ of $J(x)$ and
remark that $J(x_0)\in\mathbb R_+$ $($see Lemma \ref{lemma-Hy}$)$.

\subsection{Large deviations}\label{subsec:ldp}

In this subsection we state the LDPs for the sequence $\{\frac{n-A_n^*}{f(n)}\}_{n\in\mathbb N}$ for
different choices of the scaling function $f$. A brief summary of the results, based on the identification of three different regimes,
is given at the end of this subsection, after the statement of Theorem \ref{ldp4}.

The following theorems hold.

\begin{Theorem}\label{ldp1}
Assume \eqref{eq:hyp}, \eqref{hyp:supcrit} and set
\begin{equation}\label{fn-def}
f_1(n):=1\vee\frac{g(n)a_c^{(n)}}{n p_n},
\end{equation}
where $g$ is an arbitrary function diverging to $+\infty$, chosen
in such a way that there exists the limit $\lim_{n\to\infty}p_n f_1(n)$ $($finite or infinite$)$ and
$\lim_{n\to\infty}f_1(n)/n=\ell_1\in [0,\infty)$. Then $\{\frac{n-A_n^*}{f_1(n)}\}_{n\in\mathbb N}$ obeys an LDP
on $\overline{\mathbb R}:=\mathbb R\cup\{+\infty\}$ with speed $v_1(n):=a_c^{(n)}$ and rate function
\begin{equation*}
I_1(x):=\left\{
\begin{array}{ll}
0 &\ {\rm if}\ x=0\\
J(x_0) &\ {\rm if}\ x= \ell_1^{-1}\\
+\infty &\ {\rm otherwise},\\
\end{array}
\right .
\end{equation*}
where $\ell_1^{-1}:=+\infty$ for $\ell_1=0$.
\end{Theorem}

\begin{Remark}\label{rem:ldp1}
Let $\{X_n\}_{n\in\mathbb N}$ be a stochastic process defined by
$P(X_n=x_n^{(1)})=1-\mathfrak{p}_n$ and $P(X_n=x_n^{(2)})=\mathfrak{p}_n$, where
the sequences $\{\mathfrak{p}_n\}_{n\in\mathbb N}\subset (0,1)$ and $\{x_n^{(i)}\}_{n\in\mathbb N}\subset [0,+\infty)$, $i=1,2$, are such that
$\mathfrak{p}_n\to 0$, $x_n^{(1)}\to 0$ and $x_n^{(2)}\to\bar x\in (0,+\infty]$. Letting
$\{v_n\}_{n\in\mathbb N}\subset\R_+$ denote a sequence diverging to $+\infty$ such that  $\log\mathfrak{p}_n/v_n\to -\mathfrak c$, for some
$\mathfrak c\in (0,+\infty]$, by a direct computation, we have that  $\{X_n\}_{n\geq 1}$ obeys an LDP on $\overline{\R}$ with
speed $v_n$ and rate function
\begin{equation*}
I(x):=\left\{
\begin{array}{ll}
0 &\ {\rm if}\ x=0\\
\mathfrak c &\ {\rm if}\ x=\bar x\\
+\infty &\ {\rm otherwise}.
\end{array}
\right .
\end{equation*}
Defining
$x_n^{(1)}:=o(f_1(n))/f_1(n)$,
$x_2^{(n)}:=(n+o(n))/f_1(n)$, $\mathfrak{p}_n:=\mathrm{e}^{-J(x_0)a_c^{(n)}}$,
$\bar x:=\ell_1^{-1}$, $v_n:=a_c^{(n)}$ and $\mathfrak c:=J(x_0)$,
by Theorem \ref{ldp1} we have that
$\{(n-A_n^*)/f_1(n)\}_{n\in\mathbb N}$ obeys the same LDP as
$\{X_n\}_{n\in\mathbb N}$. As we shall discuss in Section \ref{sec:disc}, this is in accordance with an intuitive interpretation
of the result.
\end{Remark}

\begin{Theorem}\label{ldp2}
Assume \eqref{eq:hyp}, \eqref{hyp:supcrit},
\begin{equation}\label{eq:bcinfty}
b_c^{(n)}\to +\infty
\end{equation}
and let $f_2$ be a function such that $f_2(n)/b_c^{(n)}\to\ell_2\in\mathbb R_+$.
Then $\{\frac{n-A_n^*}{f_2(n)}\}_{n\in\mathbb N}$ obeys an LDP
on $\mathbb R$ with speed $v_2(n):=b_c^{(n)}$ and rate function $I_2(x):=H(\ell_2 x)$.
\end{Theorem}

\begin{Remark}\label{rem:ldp2}
Let $D_n$ denote the number of nodes in $G(n,p_n)$ with degree strictly less than $r$.
By construction, we clearly have  $n-A_n^*\ge D_n$ almost surely.
In the super-critical regime, one may naturally expect that $n-A_n^*$ behaves similarly to $D_n$, as $n\to\infty$.
At the level of the weak law of large numbers, this was pointed out by Janson's et al. in \cite{JLTV}. Indeed, if $b_c^{(n)}\to+\infty$, one may easily check that
$D_n/b_c^{(n)}\to 1$ in probability and by Theorem \ref{thm:simplifiedjanson}$(i)$ below we similarly have $(n-A_n^*)/b_c^{(n)}\to 1$ in probability.
Theorem \ref{ldp2} lifts this analogy at the level of large deviations. Indeed, under the same assumptions of Theorem \ref{ldp2},
an application of the G\"artner-Ellis theorem shows that $\{D_n/b_c^{(n)}\}_{n\in\mathbb N}$
obeys an LDP on $\R$ with speed $b_c^{(n)}$ and rate function $H$.
\end{Remark}

\begin{Theorem}\label{ldp5}
Assume \eqref{eq:hyp}, \eqref{hyp:supcrit},
\begin{equation}\label{eq:bcb}
b_c^{(n)}\to b\in (0,\infty]
\end{equation}
and let $f_3$ be a function which diverges to $+\infty$ in such a way that
\begin{equation}\label{eq:f3property}
f_3(n)/b_c^{(n)}\to+\infty\quad\text{and}\quad\frac{n p_n f_3(n)}{a_c^{(n)}}\to 0.
\end{equation}
Then $\{\frac{n-A_n^*}{f_3(n)}\}_{n\in\mathbb N}$ obeys an LDP
on $\mathbb R$ with speed $v_3(n):=-f_3(n)\log(b_c^{(n)}/f_3(n))$ and rate function
\begin{equation*}
I_3(x):=\left\{
\begin{array}{ll}
+\infty &\ {\rm if}\ x\in\R_-\\
x &\ {\rm if}\ x\geq 0.\\
\end{array}
\right .
\end{equation*}
\end{Theorem}

\begin{Theorem}\label{ldp3}
Assume \eqref{eq:hyp}, \eqref{hyp:supcrit},
\begin{equation}\label{eq:bczero}
b_c^{(n)}\to 0,
\end{equation}
\begin{equation}\label{acpninfty}
a_c^{(n)}/(n p_n)\to +\infty
\end{equation}
and let $f_4$ be a function such that
\begin{equation*}
f_4(n)\to\ell_4\in (0,\infty].
\end{equation*}

\noindent$(i)$ If $\ell_4\in\mathbb R_+$, then $\{\frac{n-A_n^*}{f_4(n)}\}_{n\in\mathbb N}$ obeys an LDP
on $\mathbb R$ with speed $v_4(n):=-\log b_c^{(n)}$ and rate function
\begin{equation}\label{J3}
I_4(x):=\left\{
\begin{array}{ll}
+\infty &\ {\rm if}\ x\in\R_-\\
\lceil\ell_4 x\rceil &\ {\rm if}\ x\geq 0.\\
\end{array}
\right .
\end{equation}

\noindent$(ii)$ If $\ell_4=+\infty$ and
\begin{equation}\label{eq:f4o}
\lim_{n\to\infty}\frac{n p_n f_4(n)}{a_c^{(n)}}=0,
\end{equation}
then $\{\frac{n-A_n^*}{f_4(n)}\}_{n\in\mathbb N}$ obeys an LDP
on $\mathbb R$ with speed $v_4(n):=-f_4(n)\log(b_c^{(n)}/f_4(n))$ and rate function $I_4:=I_3$.
\end{Theorem}

\begin{Theorem}\label{ldp4}
Assume \eqref{eq:hyp}, \eqref{hyp:supcrit},
\begin{equation}\label{eq:acpnk}
a_c^{(n)}/(n p_n)\to\gamma\in (0,\infty]
\end{equation}
and let $f_5$ be a function such that $f_5(n)\to\ell_5\in (0,\infty]$.\\
\noindent$(i)$ If $\gamma,\ell_5\in\mathbb R_+$,
then $\{\frac{n-A_n^*}{f_5(n)}\}_{n\in\mathbb N}$ obeys an LDP on
$\overline{\mathbb R}:=\mathbb R\cup\{+\infty\}$ with speed $v_5(n):=a_c^{(n)}$ and rate function
\begin{equation*}
I_5(x):=\left\{
\begin{array}{ll}
+\infty &\ {\rm if}\ x\in\R_-\\
\gamma^{-1}\lceil \ell_5 x\rceil &\ {\rm if}\ x\geq 0\\
J(x_0) &\ {\rm if}\ x=+\infty.
\end{array}
\right .
\end{equation*}
\noindent$(ii)$ If $\gamma=\ell_5=\infty$, $b_c^{(n)}\to b\in [0,\infty]$ and
$f_5(n)\sim_e\ell_5' a_c^{(n)}/(n p_n)$, for some $\ell_5'\in\mathbb R_+$, then $\{\frac{n-A_n^*}{f_5(n)}\}_{n\in\mathbb N}$ obeys an LDP on $\overline{\mathbb R}$
with speed $v_5(n):=a_c^{(n)}$ and rate function
\begin{equation*}
I_5(x):=\left\{
\begin{array}{ll}
+\infty &\ {\rm if}\ x\in\R_-\\
\ell_5'x &\ {\rm if}\ x\geq 0\\
J(x_0) &\ {\rm if}\ x=+\infty.
\end{array}
\right .
\end{equation*}
\end{Theorem}

A brief summary of the these results can be given by distinguishing the following three different regimes:
\begin{equation}\label{eq:regimes}
b_c^{(n)}\to+\infty, \qquad
b_c^{(n)}\to b\in\mathbb R_+, \qquad \text{and}
\qquad b_c^{(n)}\to 0,
\end{equation}

\begin{itemize}

\item[(1)] Under the first regime, Theorems \ref{ldp1}, \ref{ldp2}, \ref{ldp5} and \ref{ldp4} provide LDPs for
$\{\frac{n-A_n^*}{f(n)}\}_{n\in\mathbb N}$ with a divergent scaling function $f$ such that $b_c^{(n)}\lesssim f(n)\lesssim n$.
Indeed, Theorem \ref{ldp2} provides an LDP with a scaling function $f(n)=f_2(n)\sim b_c^{(n)}$;
Theorem \ref{ldp5} provides an LDP with a scaling function $f=f_3$ such that $b_c^{(n)}\ll f_3(n)\ll a_c^{(n)}/(n p_n)$;
Theorem \ref{ldp4}$(ii)$ provides an LDP with a scaling function $f(n)=f_5(n)\sim a_c^{(n)}/(n p_n)$;
Theorem \ref{ldp1} provides an LDP with a scaling function $f=f_1$ such that
$a_c^{(n)}/(n p_n)\ll f_1(n)\lesssim n$.
\item[(2)] Under the second regime, Theorems \ref{ldp1}, \ref{ldp5} and \ref{ldp4} provide LDPs for
$\{\frac{n-A_n^*}{f(n)}\}_{n\in\mathbb N}$ with a
divergent scaling function $f$ such that $f(n)\lesssim n$.
Indeed, Theorem \ref{ldp5} provides an LDP with a scaling function $f=f_3$ such that
$f_3(n)\ll a_c^{(n)}/(n p_n)$; Theorem \ref{ldp4}$(ii)$ provides an LDP with a scaling function $f(n)=f_5(n)\sim a_c^{(n)}/(n p_n)$;
Theorem \ref{ldp1} provides an LDP with a scaling function $f=f_1$  such that
$a_c^{(n)}/(n p_n)\ll f_1(n)\lesssim n$.
\item[(3)] Under the third regime, Theorems \ref{ldp1}, \ref{ldp3} and \ref{ldp4} provide LDPs for
$\{\frac{n-A_n^*}{f(n)}\}_{n\in\mathbb N}$ with a scaling function which may be
either convergent or divergent. More  precisely, we  distinguish the following three cases:
\begin{equation}\label{eq:cases}
a_c^{(n)}/(np_n)\to+\infty, \qquad  a_c^{(n)}/(np_n)\to\gamma\in\mathbb R_+ \qquad \text{and} \qquad  a_c^{(n)}/(np_n)\to 0.
\end{equation}

In the first case, Theorem \ref{ldp3}$(i)$ provides an LDP with a convergent scaling function $f(n)=f_4(n)\sim_e\ell_4$ and
Theorems \ref{ldp3}$(ii)$, \ref{ldp4}$(ii)$ and \ref{ldp1} provide LDPs with a divergent scaling function $f$ such that
$f(n)\lesssim n$. Indeed, Theorems \ref{ldp3}$(ii)$
provides an LDP with a divergent scaling function $f=f_4$ such that $f_4(n)\ll a_c^{(n)}/(n p_n)$;
Theorem \ref{ldp4}$(ii)$ provides an LDP with a divergent scaling function $f(n)=f_5(n)\sim a_c^{(n)}/(n p_n)$;
Theorem \ref{ldp1} provides an LDP with a divergent scaling function $f=f_1$ such that
$a_c^{(n)}/(n p_n)\ll f_1(n)\lesssim n$.

In the second case, Theorem \ref{ldp4}$(i)$ provides an LDP with a
convergent scaling function $f(n)=f_5(n)\sim_e\ell_5$ and Theorem \ref{ldp1}  provides an LDP with a divergent scaling function $f=f_1$
such that $f_1(n)\lesssim n$.

In the third case, Theorem \ref{ldp1} covers the whole range providing LDPs
with a scaling function $f=f_1$ which is either convergent, i.e. $f_1(n)\sim_e\ell$, for some constant $\ell\geq 1$, or divergent
in such a way that $f_1(n)\lesssim n$.
\end{itemize}

We conclude this subsection recalling (for later purposes) that, as noticed in \cite{JLTV} (see formula 3.10 therein), under \eqref{eq:hyp},
the three regimes \eqref{eq:regimes} are equivalent to
\begin{equation}\label{rfIexpression}
n p_n-(\log n +(r-1)\log\log n)\to\left\{
\begin{array}{ll}
-\infty&\\
-\log((r-1)!b)\\
+\infty,&\\
\end{array}
\right .
\end{equation}
respectively.

\subsection{Laws of large numbers}

The following laws of large numbers (LLNs) are corollaries of the previous LDPs.
Their proofs are omitted since they are based on a standard application of the large deviation estimates
and the Borel-Cantelli lemma.

\begin{Theorem}\label{lln2}
Under the assumptions and notation of Theorem \ref{ldp2}, we have
\[
\frac{n-A_n^*}{f_2(n)}\to\ell_2^{-1},\quad\text{almost surely.}
\]
provided that $v_2(n)\gg log n$.
\end{Theorem}

\begin{Theorem}\label{lln5}
Under the assumptions and notation of Theorem \ref{ldp5} $($Theorem \ref{ldp3}, Theorem \ref{ldp4}, respectively$)$, we have
\begin{equation}\label{eq:divg}
\frac{n-A_n^*}{f_3(n)}\to 0,\quad\text{almost surely}
\end{equation}
$($$(n-A_n^*)/f_4(n)\to 0$, $(n-A_n^*)/f_5(n)\to 0$, almost surely, respectively$)$  provided that  $v_3(n)\gg \log n$. ($v_4(n)\gg \log n $,  $v_5(n)\gg \log n$). 
\end{Theorem}

To better position Theorems \ref{lln2} and \ref{lln5} with respect to the corresponding results in \cite{JLTV}, in Theorem \ref{thm:simplifiedjanson} below we state
the main achievements of \cite{JLTV} for the bootstrap percolation process in the super-critical regime.

\begin{Theorem}\label{thm:simplifiedjanson}
Assume \eqref{eq:hyp} and \eqref{hyp:supcrit}. Then:\\
$(i)$ If
$b_c^{(n)}\to+\infty$, then $(n-A_n^*)/b_c^{(n)}\to 1$, in probability.\\
$(ii)$ If
$b_c^{(n)}\to 0$, then $P(A_n^*=n)\to 1$.\\
$(iii)$ If
$b_c^{(n)} \to b\in\mathbb R_+$, then $n-A_n^*$ converges in distribution to a Poisson distributed random variable with mean $b$.
\end{Theorem}

Taking $f_2(n)=b_c^{(n)}$  (under the assumption that $b_c^{(n)}\gg \log n $) in Theorem \ref{lln2}, we have
$(n-A_n^*)/b_c^{(n)}\to 1$, almost surely. This LLN generalizes the one in
Theorem \ref{thm:simplifiedjanson}$(i)$, where
the convergence holds in probability.

Taking either $f_5\equiv 1$ in Theorem \ref{lln5}, we have $n-A_n^*\to 0$, almost surely.
This LLN generalizes the one in Theorem \ref{thm:simplifiedjanson}$(ii)$. Indeed, since the random variable $n-A_n^*$
takes values in $\mathbb N\cup\{0\}$, for any $\varepsilon\in (0,1)$,
\[
P(A_n^*\neq n)=P(n-A_n^*>0)=P(n-A_n^*>\varepsilon),\quad n\in\mathbb N.
\]
Therefore, if $n-A_n^*\to 0$ almost surely, then $n-A_n^*\to 0$ in probability, and so
$P(A_n^*=n)\to 1$.

If $b_c^{(n)}\to b\in\mathbb R_+$, then by Theorem \ref{thm:simplifiedjanson}$(iii)$, we have that $n-A_n^*\to\mathrm{Po}(b)$ in distribution, where
$\mathrm{Po}(b)$ denotes a Poisson distributed random variable with mean $b$. It follows that, for any divergent function $f_3$, we have
$(n-A_n^*)/f_3(n)\to 0$, in probability. This weak LLN is extended by formula \eqref{eq:divg} in Theorem \ref{lln5}.

\subsection{Tail asymptotics}\label{subsec:tail}

The proofs of Theorems \ref{ldp1}, \ref{ldp2}, \ref{ldp5}, \ref{ldp3} and \ref{ldp4} are based, respectively,
on Propositions \ref{ldp:supcrit}, \ref{ldp:supcritbis}, \ref{ldp:supcritpenta}, \ref{ldp:supcritter} and \ref{ldp:supcritquater} below,
which provide asymptotic estimates for the tail of $n-A_n^*$. For reader's convenience, we summarize these tail estimates
in Tables \ref{tab:1}-\ref{tab:5}, reported after the statement of Proposition \ref{ldp:supcritquater}. An informal discussion of these results
is postponed to Section \ref{sec:disc}.

\begin{Proposition}\label{ldp:supcrit}
Under the same notation and assumptions of Theorem \ref{ldp1}, for an arbitrarily chosen  $\varepsilon\in (0,\ell_1^{-1})$ we have:
\begin{equation*}
\lim_{n\to\infty}\frac{1}{a_c^{(n)}}\log P\left(\frac{n-A_n^*}{f_1(n)}>\varepsilon\right)=-J(x_0).
\end{equation*}
\end{Proposition}

\begin{Proposition}\label{ldp:supcritbis},
Under the same notation and assumptions of Theorem \ref{ldp2}, for any arbitrarily chosen $\varepsilon>\ell_2^{-1}$, we have:
\begin{equation}\label{eq:lim1}
\lim_{n\to\infty}\frac{1}{b_c^{(n)}}\log P\left(\frac{n-A_n^*}{f_2(n)}>\varepsilon\right)=-H(\ell_2\varepsilon)
\end{equation}
and, for any arbitrarily fixed $\varepsilon\in (0,\ell_2^{-1})$,
\begin{equation}\label{eq:lim2}
\lim_{n\to\infty}\frac{1}{b_c^{(n)}}\log P\left(\frac{n-A_n^*}{f_2(n)}\leq\varepsilon\right)=-H(\ell_2\varepsilon).
\end{equation}
\end{Proposition}

\begin{Proposition}\label{ldp:supcritpenta}
Under the same notation and assumptions of Theorem \ref{ldp5}, for an arbitrarily  chosen  $\varepsilon\in\mathbb R_+$, we have:
\begin{equation*}
\lim_{n\to\infty}\frac{1}{-f_3(n)\log(b_c^{(n)}/f_3(n))}\log P\left(\frac{n-A_n^*}{f_3(n)}>\varepsilon\right)=-\varepsilon.
\end{equation*}
\end{Proposition}

\begin{Proposition}\label{ldp:supcritter}
Under the same notation and assumptions of Theorem \ref{ldp3}, for an arbitrarily chosen $\varepsilon\in\mathbb R_+$, we have:\\
\noindent$(i)$ if $\ell_4\in\mathbb R_+$, then
\begin{equation}\label{eq:lim1ter}
\lim_{n\to\infty}\frac{1}{-\log b_c^{(n)}}\log P\left(\frac{n-A_n^*}{f_4(n)}>\varepsilon\right)=-\lceil\ell_4\varepsilon\rceil;
\end{equation}
\noindent$(ii)$ if $\ell_4=\infty$ and \eqref{eq:f4o} holds, then
\begin{equation}\label{eq:lim1terbis}
\lim_{n\to\infty}\frac{1}{-f_4(n)\log (b_c^{(n)}/f_4(n))}\log P\left(\frac{n-A_n^*}{f_4(n)}>\varepsilon\right)=-\varepsilon.
\end{equation}
\end{Proposition}

\begin{Proposition}\label{ldp:supcritquater}
Under the same notation and assumptions of Theorem \ref{ldp4}, for an arbitrarily chosen $\varepsilon\in\mathbb R_+$, we have:\\
\noindent$(i)$ if $\gamma,\ell_5\in\mathbb R_+$, then
\begin{equation*}
\lim_{n\to\infty}\frac{1}{a_c^{(n)}}\log P\left(\frac{n-A_n^*}{f_5(n)}>\varepsilon\right)=-J(x_0)\wedge\gamma^{-1}\lceil\ell_5\varepsilon\rceil;
\end{equation*}
\noindent$(ii)$ if $\gamma=\ell_5=+\infty$, $b_c^{(n)}\to b\in [0,\infty]$ and $f_5(n)\sim_e\ell_5'a_c^{(n)}/(n p_n)$, for some $\ell_5'\in\mathbb R_+$,
then
\begin{equation*}
\lim_{n\to\infty}\frac{1}{a_c^{(n)}}\log P\left(\frac{n-A_n^*}{f_5(n)}>\varepsilon\right)=-J(x_0)\wedge \ell_5'\varepsilon.
\end{equation*}
\end{Proposition}

Using an obvious notation, the asymptotic estimates for the right tail of
$n-A_n^*$, provided by the above propositions,
can be summarized as
\begin{equation}\label{eq:genasrel}
\log P\left(\frac{n- A_n^*}{f(n)}>\varepsilon\right)\sim_e-\mathcal{I}(\varepsilon)v(n),\qquad\text{$\varepsilon>z_0$}
\end{equation}
where
$z_0$ is the most probable value of $\{\frac{n- A_n^*}{f(n)}\}_{n\in\mathbb N}$ as $n\to \infty$.
Tables \ref{tab:1}-\ref{tab:5} below report the function $\mathcal I$
and the speed function
$v$ for different choices of the scaling function $f$, in the three different regimes \eqref{eq:regimes}.

\begin{table}[h]
\begin{center}
\vspace{0mm}
\caption{The functions $\mathcal{I}$, $v$ and $f$,
in the regime $b^{(n)}_c \to +\infty$\label{tab:1}}{
\begin{tabular}{|c||c|c|c|c|}
\hline
\Top
$b^{(n)}_c \to +\infty$
& \multicolumn{4}{|c|}{divergent $f(n)$}
\Bottom \\
\cline{2-5}
&
\Top
$f(n) \sim_e \ell b^{(n)}_c$ & $b^{(n)}_c \ll f(n) \ll \frac{a_c^{(n)}}{n p_n}$ &
$f(n) \sim_e \ell\frac{a_c^{(n)}}{n p_n}$ & $\frac{a_c^{(n)}}{n p_n}\ll f(n)\lesssim n$
\Bottom \\
\hline
\hline
\Top
 $v(n)$ & $b^{(n)}_c$ & $-f(n)\log(b_c^{(n)}/f(n))$ & $a_c^{(n)}$ & $a_c^{(n)}$
\Bottom \\
\hline
\Top
$\mathcal I(\varepsilon)$ & $H(\ell\varepsilon)$ & $\varepsilon$ & $J(x_0)\wedge\ell\varepsilon$ & $J(x_0)$
\Bottom \\
\hline
Prop. & \ref{ldp:supcritbis} & \ref{ldp:supcritpenta} & \ref{ldp:supcritquater} & \ref{ldp:supcrit} \\
\hline
\end{tabular}}
\end{center}
\vspace{-0mm}
\end{table}

\begin{table}
\begin{center}
\vspace{0mm}
\caption{The functions $\mathcal I$, $v$ and $f$, in the regime $b^{(n)}_c \to b\in\mathbb R_+$ \label{tab:2}}{%
\begin{tabular}{|c||c|c|c|}
\hline
\Top
\multirow{2}{*}{$b^{(n)}_c \to b\in\mathbb R_+$} & \multicolumn{3}{|c|}{divergent $f(n)$}
\Bottom \\
\cline{2-4}
&  \Top $f(n) \ll \frac{a_c^{(n)}}{n p_n}$ &
$f(n) \sim_e \ell\frac{a_c^{(n)}}{n p_n}$ & $\frac{a_c^{(n)}}{n p_n} \ll f(n)\lesssim n$
\Bottom \\
\hline
\hline
\Top
$v(n)$ & $f(n)\log f(n)$ & $a_c^{(n)}$ & $a_c^{(n)}$
\Bottom \\
\hline
\Top
 $\mathcal I(\varepsilon)$ & $\varepsilon$ & $J(x_0)\wedge\ell\varepsilon$ & $J(x_0)$
\Bottom \\
\hline
Prop. & \ref{ldp:supcritpenta} & \ref{ldp:supcritquater} & \ref{ldp:supcrit} \\
\hline
\end{tabular}}
\end{center}
\vspace{-0mm}
\end{table}

\begin{table}
\begin{center}
\vspace{0mm}
\caption{The functions $\mathcal{I}$, $v$ and $f$,
in the regime $b^{(n)}_c \to 0$, $\frac{a_c^{(n)}}{n p_n} \to +\infty$ \label{tab:3}}{
\begin{tabular}{|c||c|c|c|c|}
\hline
\Top
$b^{(n)}_c \to 0$ & \multicolumn{4}{|c|}{$f(n)$}
\Bottom \\
\cline{2-5}
$\frac{a_c^{(n)}}{n p_n} \to +\infty$ &
\Top
$f(n) \sim_e \ell$ & $f$ divergent: $f(n) \ll \frac{a_c^{(n)}}{n p_n}$ &
$f(n) \sim_e \ell\frac{a_c^{(n)}}{n p_n}$ & $\frac{a_c^{(n)}}{n p_n}\ll f(n)\lesssim n$
\Bottom \\
\hline
\hline
\Top
 $v(n)$ & $-\log b_c^{(n)}$ & $-f(n)\log(b_c^{(n)}/f(n))$ & $a_c^{(n)}$ & $a_c^{(n)}$
\Bottom \\
\hline
\Top
$\mathcal I(\varepsilon)$ & $\lceil\ell\varepsilon\rceil$ & $\varepsilon$ & $J(x_0)\wedge\ell\varepsilon$ & $J(x_0)$
\Bottom \\
\hline
Prop. & \ref{ldp:supcritter} & \ref{ldp:supcritter} & \ref{ldp:supcritquater} & \ref{ldp:supcrit} \\
\hline
\end{tabular}}
\end{center}
\vspace{-0mm}
\end{table}

\begin{table}
\begin{center}
\vspace{0mm}
\caption{The functions $\mathcal{I}$, $v$ and $f$,
in the regime $b^{(n)}_c \to 0$, $\frac{a_c^{(n)}}{n p_n} \to \gamma\in\mathbb R_+$ \label{tab:4}}{%
\begin{tabular}{|c||c|c|}
\hline
\Top
$b^{(n)}_c \to 0$ & \multicolumn{2}{|c|}{$f(n)$}
\Bottom \\
\cline{2-3}
$\frac{a_c^{(n)}}{n p_n} \to \gamma\in\mathbb R_+$ & \Top $f(n) \sim_e \ell$ &  $f$ divergent: $f(n)\lesssim n$
\Bottom \\
\hline
\hline
\Top
$v(n)$ & $a_c^{(n)}$ & $a_c^{(n)}$
\Bottom \\
\hline
\Top
$\mathcal{I}(\varepsilon)$ & $J(x_0)\wedge\gamma^{-1}\lceil\ell\varepsilon\rceil$ & $J(x_0)$
\Bottom \\
\hline
Prop. & \ref{ldp:supcritquater} & \ref{ldp:supcrit} \\
\hline
\end{tabular}}
\end{center}
\vspace{-0mm}
\end{table}

\begin{table}
\begin{center}
\vspace{0mm}
\caption{The functions $\mathcal{I}$, $v$ and $f$, in the regime $b^{(n)}_c \to 0$, $\frac{a_c^{(n)}}{n p_n} \to 0$ \label{tab:5}}{%
\begin{tabular}{|c||c|}
\hline
\Top
$b^{(n)}_c \to 0$ & \multicolumn{1}{|c|}{$f(n)$}
\Bottom \\
\cline{2-2}
$\frac{a_c^{(n)}}{n p_n} \to 0$ & \Top $f(n)\sim_e\ell\geq 1$ or $f$ divergent: $f(n)\lesssim n$
\Bottom \\
\hline
\hline
\Top
$v(n)$ & $a_c^{(n)}$
\Bottom \\
\hline
\Top
$\mathcal{I}(\varepsilon)$ & $J(x_0)$
\Bottom \\
\hline
Prop. & \ref{ldp:supcrit} \\
\hline
\end{tabular}}
\end{center}
\vspace{-0mm}
\end{table}

\newpage

\section{Informal discussion}\label{sec:disc}

In this section we provide an informal explanation of our results.

{Formula  \eqref{eq:genasrel} establishes an asymptotic relationship
between the \lq\lq resolution'' $f$
at which we observe the fluctuations of the random
variable $n-A_n^*$
and the decay rate $v$ at which the associated tail probability vanishes.
One might intuitively expect the following \lq\lq monotonicity'' property:
given two scaling functions  $\widetilde{f}$ and  $f$, with speed functions
$\widetilde{v}$ and $v$, respectively, if
$\widetilde{f}(n)\ll f(n)$ then $\widetilde{v}(n)\ll v(n)$, i.e., to a larger \lq\lq resolution''
corresponds a faster decay rate.

Propositions \ref{ldp:supcrit}, \ref{ldp:supcritbis}, \ref{ldp:supcritpenta}, \ref{ldp:supcritter} and \ref{ldp:supcritquater}
reveal that the \lq\lq monotonicity'' property suggested by the intuition is correct only as long as the scaling function $f$
is such that $f(n)\ll a_c^{(n)}/(n p_n)$ and $a_c^{(n)}/(n p_n)\to+\infty$. Specifically, Propositions \ref{ldp:supcritbis}, \ref{ldp:supcritpenta} and \ref{ldp:supcritter} support
the intuition but Propositions \ref{ldp:supcrit} and \ref{ldp:supcritquater} provide a somehow  counterintuitive result
showing that if the scaling function $f$ is chosen in such a way that $a_c^{(n)}/(n p_n)\ll f(n)$ and $a_c^{(n)}/(n p_n)\to+\infty$,
then the speed $v$ becomes insensitive to the selected scaling function, indeed $v(n)=a_c^{(n)}$.
In conclusion, our results reveal that the bootstrap percolation process exhibits a fairly different behavior according to
either $f(n)\ll a_c^{(n)}/(n p_n)$ or $a_c^{(n)}/(n p_n)\ll f(n)$.

To intuitively explain the reason of such behavior,
we note that the average number of usable nodes
evolves, as the time $t$ increases, according to the function
$e_n(t):=E[A_n(t)]-t$, $a_n\leq t\leq n$. 
As shown in ~\cite{JLTV}, starting at $t=a_n$,
$e_n(\cdot)$ first decreases up to reach a minimal value $\simeq a_n -a_c^{(n)}$ at $t\simeq a_c^{(n)}$; then $e_n(\cdot)$ increases quickly (super-linearly)
up to reach a value $\simeq n-t$, where $t$ is such that $E[A_n(t)]\simeq n$; finally, $e_n(\cdot)$ decreases linearly
and approaches zero at $t\simeq n-b_c^{(n)}\sim_e n$.

Since the bootstrap percolation process stops the first time the number of active and not yet used nodes equals zero,
assuming that this quantity is sufficiently concentrated around its average we expect the bootstrap percolation
process to stop either at a critical time $t\simeq a_c^{(n)}\sim a_n$ or when the process gets sufficiently close to its natural termination, i.e., at a critical time
$t\simeq n- b_c^{(n)}\sim_e n$. Since we are considering $T_n=A_n^*$ only on time intervals of the form
$[a_n,n-\varepsilon f(n))$, the choice of $f(n)$ has a direct impact on the probability that the process stops
before $n-\varepsilon f(n)$ at a time $t=n-o(n)$, but it has no effect on the probability that the process stops at a time  $t\simeq a_c^{(n)}\sim a_n$.

It turns out that if $f$ is such that $ \frac{a_c^{(n)}}{n p_n}\ll f(n)$ (recall that,
under the assumptions \eqref{eq:hyp}, $b_c^{(n)}\ll a_c^{(n)}/(n p_n)$), then the  probability that the bootstrap percolation process stops at
an early stage (i.e., at a time $t\simeq a_c^{(n)}$) is dominating.
Instead, if we choose $f$ so that $f(n)\ll\frac{a_c^{(n)}}{n p_n}$,
then probability that the process stops just before $n-\varepsilon f(n)$ (i.e.,
at a time $t\simeq n-f(n)$) becomes dominating. Finally, if $f(n)\sim a_c^{(n)}/(n p_n)$,
then the probability that the process stops just before $n-\varepsilon f(n)$  becomes
comparable with the probability that the process stops at a time $t\simeq a_c^{(n)}$.

More precisely, our analysis shows that by selecting a scaling function $f$ such that  $a_c^{(n)}/(n p_n)\ll f(n)$, we have
\[
\log P(T_n<n-\varepsilon f(n))\sim_e \log P(\exists t\in\{a_n,\ldots,\lfloor K a_c^{(n)}\rfloor\}:\,\,A_n(t)\leq t)\sim_e
- J(x_0)a_c^{(n)}
\]
and
\[
\limsup_{n\to \infty}\frac{1}{a_c^{(n)}}\log P(\exists t\in [n-\varepsilon' f(n),n-\varepsilon f(n)]:\,\,A_n(t)\leq t)
=-\infty,
\]
for any $\varepsilon'>\varepsilon$ and any $K\in\mathbb R_+$ big enough (see the proof of Proposition \ref{ldp:supcrit})).

Instead, if we choose $f$ in such a way that $f(n)\ll a_c^{(n)}/(n p_n)$, then
there exist two functions $v$ and $\mathcal I$ with $v(n)\ll a_c^{(n)}$ such that
\[
\log P(T_n<n-\varepsilon f(n))\sim_e \log P(\exists t\in\{\lfloor n- \varepsilon' f(n)\rfloor,\ldots,\lfloor n-\varepsilon f(n)\rfloor\}:\,\,A_n(t)\leq t)\sim_e
-\mathcal{I}(\varepsilon)v(n)
\]
and
\[
\lim_{n\to\infty}\frac{1}{v(n)}\log P\left(\exists t\in\{a_n,\ldots,\lfloor K a_c^{(n)}\rfloor\}:\,\,A_n(t)\leq t\right)=-\infty,
\]
for any $\varepsilon'>\varepsilon$ and any
$K\in\mathbb R_+$  (see the proofs of  Propositions \ref{ldp:supcritbis}, \ref{ldp:supcritpenta}, \ref{ldp:supcritter}).

Finally, if $f(n)\sim a_c^{(n)}/(n p_n)$, then
\[
\log P(T_n<n-\varepsilon f(n))\sim_e\log P(\exists t\in\{\lfloor n-\varepsilon'\ f(n)\rfloor,\ldots,\lfloor n-\varepsilon f(n)\rfloor\}:\,\,A_n(t)\leq t)\sim
-\varepsilon a_c^{(n)},
\]
for any $\varepsilon'>\varepsilon$, and again
\[
\log P(\exists t\in\{a_n,\ldots,\lfloor K a_c^{(n)}\rfloor\}:\,\,A_n(t)\leq t)\sim_e
- J(x_0)a_c^{(n)}
\]
(see the proof of Proposition  \ref{ldp:supcritquater}).}

\section{Proofs}\label{sec:proofs}

\subsection{Asymptotic relations and deviation bounds}

We start stating some preliminary asymptotic relations and deviation bounds, that will come in handy in the proofs.

\subsubsection{Asymptotic relations concerning the binomial distribution}

Let $\{q_n\}_{n\in\mathbb N}\subset (0,1)$, $\{m_n\}_{n\in\mathbb N}\subset [1,+\infty)$, $\{r_n\}_{n\in\mathbb N}\subset\mathbb N$ and suppose $q_n\to 0$. The
following asymptotic relations hold.

If $q_n m_n\to 0$ and $m_n\to\infty$, then, for any $k\in\mathbb N$,
\begin{equation}\label{eq:Jasonbin}
P(\mathrm{Bin}(\lfloor m_n\rfloor,q_n)\geq k)=\frac{(q_n m_n)^{k}}{k!}(1+o(1))
\end{equation}
(see e.g. formula $(8.1)$ in \cite{JLTV}).

If $r_n\to\infty$, the limit $\lim_{n\to\infty}q_n m_n$ exists (finite or not), $\frac{r_n}{q_n m_n}\to +\infty$ and $r_n m_n^{-1}\to 0$, then
\begin{equation}\label{eq:nuovorate}
\log P(\mathrm{Bin}(\lfloor m_n\rfloor,q_n)\geq r_n)=r_n\log\left(\frac{m_n q_n}{r_n}\right)(1+o(1)).
\end{equation}
We refer the reader to the Appendix for the proof of \eqref{eq:Jasonbin} and \eqref{eq:nuovorate}.

If $q_n m_n\to\infty$, then, for an arbitrarily fixed $k\in\mathbb N$,
\begin{equation}\label{eq:svilpi}
1-P(\mathrm{Bin}(\lfloor m_n\rfloor,q_n)\geq k)=(1-q_n)^{m_n}\frac{(q_n m_n)^{k-1}}{(k-1)!}(1+o(1))
\end{equation}
(see e.g. formula $(3.7)$ in \cite{JLTV}).

\subsubsection{Asymptotic relations concerning $a_c^{(n)}$, $\pi_n(t)$ and $b_c^{(n)}$}

One may easily verify that
\begin{align}
a_c^{(n)}\sim\frac{1}{(p_n n^{1/r})^{r/(r-1)}}.\label{eq:tc}
\end{align}
By the definition of $\pi_n$, \eqref{eq:hyp}, \eqref{eq:Jasonbin} and the definition of $a_c^{(n)}$, for any fixed $x\in\R_+$, we have
\begin{align}
n\pi_n(x a_c^{(n)})&\sim_e n\frac{x^{r}}{r!}(p_n a_c^{(n)})^{r}\nonumber\\
&=x^{r}(a_c^{(n)})^{r-1}\frac{n p_n^{r}}{r!}a_c^{(n)}\nonumber\\
&=x^{r}\left(1-\frac{1}{r}\right)^{r-1}\frac{(r-1)!}{np_n^{r}}\frac{np_n^{r}}{r!}a_c^{(n)}\nonumber\\
&=\frac{1}{r}\left(1-\frac{1}{r}\right)^{r-1}x^{r}a_c^{(n)}.\label{eq:speed}
\end{align}
Moreover, under \eqref{eq:hyp}, for any $f$ such that $f(n)=o(p_n^{-1})$, by \eqref{eq:svilpi}
\begin{equation}\label{eq:nunomenopi}
1-\pi_n(n-f(n))\sim_e\frac{b_c^{(n)'}}{n}\to 0,
\end{equation}
where
\[
b_c^{(n)'}:=n\frac{(n p_n)^{r-1}}{(r-1)!}(1-p_n)^{n}.
\]
For later purposes, we remark that if $\lim_{n\to\infty}b_c^{(n)}=b\in (0,\infty]$, then
\begin{equation}\label{eq:logblogbc1}
b_c^{(n)}=b_c^{(n)'}(1+o(1)),
\end{equation}
and that if $\lim_{n\to\infty}b_c^{(n)}=0$, then
\begin{equation}\label{eq:logblogbc}
\log b_c^{(n)}\sim_e \log b_c^{(n)'}\quad\text{and}\quad\log b_c^{(n)}\sim -n p_n,
\end{equation}
see the Appendix for a proof of \eqref{eq:logblogbc1} and \eqref{eq:logblogbc}.

\subsubsection{Deviation bounds for the binomial distribution}\label{penrose}
Throughout this paper we will extensively exploit some
classical deviation bounds for the binomial distribution (see e.g. Lemma 1.1 p. 16 in \cite{P}), which we report here
for the sake of completeness.

Let the function $H$ be defined by \eqref{eq:Hx} and set $\mu:=n p$, for $n\in\mathbb{N}$, $p\in(0,1)$. For any
$0<k<n$, we have:\\
\noindent if $k\geq\mu$, then
\begin{equation}\label{Penrose1.5}
P(\mathrm{Bin}(n,p)\geq k)\leq\exp\left(-\mu H\left(\frac{k}{\mu}\right)\right);
\end{equation}
if $k\leq\mu$, then
\begin{equation}\label{Penrose1.6}
P(\mathrm{Bin}(n,p)\leq k)\leq\exp\left(-\mu H\left(\frac{k}{\mu}\right)\right);
\end{equation}
if $k\geq\mathrm{e}^2\mu$, then
\begin{equation}\label{Penrose1.7}
P(\mathrm{Bin}(n,p)\geq k)\leq\exp\left(-\left(\frac{k}{2}\right)\log\left(\frac{k}{\mu}\right)\right).
\end{equation}

\subsection{Proofs of Proposition \ref{ldp:supcrit} and Theorem \ref{ldp1}}

The proof of Theorem \ref{ldp1} is based on Proposition \ref{ldp:supcrit}, whose proof exploits in turn the following lemmas:

\begin{Lemma}\label{ldp}
Assume $a_n=o(n)$ and let $\{\pi(n)\}_{n\in\mathbb N}\subset (0,1)$ and $\{v(n)\}_{n\in\mathbb N}\subset\mathbb R_+$
be two sequences such that $v(n)=o(n)$ and $n\pi(n)\sim_e v(n)$. Then
$\biggl\{\frac{\mathrm{Bin}(n-a_n,\pi(n))}{v(n)}\biggr\}_{n\in\mathbb N}$ obeys an LDP on $\mathbb R$
with speed $v(n)$ and rate function $I:=H$.
\end{Lemma}

\begin{Lemma}\label{lemma-Hy}
We have$:$\\
\noindent$(i)$ $x<h(x)$ for any $x\geq 0$ $($and so $J$ is strictly positive on $[0,\infty)$$)$, whenever $\alpha>1$.\\
\noindent$(ii)$ $J$ admits a unique point of minimum $x_0\in\R_+$ on $[0,\infty)$.\\
\end{Lemma}

The proofs of these lemmas are given in the latter paragraph of this subsection.

\subsubsection{Proof of Proposition \ref{ldp:supcrit}}
As a guide to the intuition, we start by briefly describing the outline of the proof.
For any $\varepsilon\in (0,\ell_1^{-1})$ and $n\in\mathbb N$
large enough, we have from \eqref{eq:T} that :
\begin{align}
P\left(\frac{n-A_n^*}{f_1(n)}\leq\varepsilon\right)=P(A_n^*\geq n-\varepsilon f_1(n))=P(A_n(t)>t,\quad\forall t=a_n,\ldots,\lfloor n-\varepsilon f_1(n)\rfloor).\nonumber
\end{align}
Therefore,
\begin{align}
P\left(\frac{n-A_n^*}{f_1(n)}>\varepsilon\right)&=1-P(A_n(t)>t,\quad\forall t=a_n,\ldots,\lfloor n-\varepsilon f_1(n)\rfloor)\nonumber\\
&=P\left(\bigcup_{t=a_n}^{\lfloor n-\varepsilon f_1(n)\rfloor}\{A_n(t)-t\leq 0\}\right)\nonumber\\
&=P\left(\bigcup_{t=a_n}^{\lfloor n-\varepsilon f_1(n)\rfloor}\{S_n(t)+a_n-t\leq 0\}\right),\label{eq:pdopo}
\end{align}
where we used \eqref{eq:A(t)}. We rewrite the event
$\{(n-A_n^*)/f_1(n)>\varepsilon\}$ as the union of the events
\begin{equation}\label{eq:nuovoB1}
B_1^{(n)}:=\bigcup_{t=a_n}^{\lfloor K a_c^{(n)}\rfloor}\{S_n(t)+a_n-t\leq 0\},
\end{equation}
\begin{equation}\label{eq:nuovoB2}
B_2^{(n)}:=\bigcup_{t=\lceil K a_c^{(n)}\rceil}^{\lfloor p_n^{-1}\rfloor}\{S_n(t)+a_n-t\leq 0\},
\end{equation}
\begin{equation}\label{eq:b3}
B_3^{(n)}:=\bigcup_{t=\lfloor p_n^{-1}\rfloor}^{\lfloor n-\varepsilon f_1(n)\rfloor}\{S_n(t)+a_n-t\leq 0\},
\end{equation}
where, as it will be made precise later on, $K$ is a sufficiently large constant. Note that by construction
$B_1^{(n)}=\{ T_n\le K a_c^{(n)}\}=\{ A^*_n\le K a_c^{(n)}\}$, while $B_2^{(n)}\supseteq \{ K a_c^{(n)} \le T_n\le\lfloor p_n^{-1} \rfloor\}=\{ K a_c^{(n)} \le A_n^*\le\lfloor p_n^{-1} \rfloor\}$
and $B_3^{(n)}\supseteq \{ \lfloor p_n^{-1}\rfloor \le T_n\le\lfloor n- \varepsilon f_1(n) \rfloor\} =  \{ \lfloor p_n^{-1}\rfloor \le A^*_n\le\lfloor n- \varepsilon f_1(n) \rfloor\}$.

Basically  in the proof we show  that $B_1^{(n)}$
is the dominating event and  we provide  tight asymptotic estimates  for  $P(B_1^{(n)})$.
More precisely, since, by construction:
\begin{equation}\label{eq:inequalities}
P(B_1^{(n)})\leq P\left(\frac{n-A_n^*}{f_1(n)}>\varepsilon\right)\leq P(B_1^{(n)})+P(B_2^{(n)})+P(B_3^{(n)}),
\end{equation}
the claim will follow by the principle of the largest term (see e.g. Lemma 1.2.15 p. 7 in \cite{DZ}) provided that we are able to show that:
\begin{equation}\label{eq:compasB1}
\lim_{n\to\infty}\frac{1}{a_c^{(n)}}\log P(B_1^{(n)})=-J(x_0),
\end{equation}
\begin{equation}\label{eq:compasb2}
\limsup_{n\to\infty}\frac{1}{a_c^{(n)}}\log P(B_2^{(n)})\leq-C_K,\quad\text{for some constant $C_K>J(x_0)$}
\end{equation}
and
\begin{equation}\label{eq:compasb3}
\lim_{n\to\infty}\frac{1}{a_c^{(n)}}\log P(B_3^{(n)})=-\infty.
\end{equation}
The proofs of \eqref{eq:compasB1}, \eqref{eq:compasb2} and \eqref{eq:compasb3} are based
on the binomial structure of $S_n(t)$ and $\pi_n(t)$, which allows to exploit
Lemma \ref{ldp} and
the deviation bounds summarized in Section \ref{penrose}.

{We proceed by dividing the proof in four steps.
 In the first step, starting from the LDP principle stated in Lemma \ref{ldp}, we derive a new LDP for the  sequence
$\{S_n(\kappa_n(x))/((1-r^{-1})^{-1}h(x)a_c^{(n)})\}_{n\in\mathbb N}$,  where:
\[
\kappa_n(x):=(\alpha+(1-r^{-1})^{-1}x+o(1))a_c^{(n)},\quad\text{$x\geq 0$, $n\in\mathbb N$.}
\]
In the second step,  we  employ the  previously obtained LDP  to prove \eqref{eq:compasB1}.
In the third step, we prove \eqref{eq:compasb2}. At last, in the fourth step we  prove \eqref{eq:compasb3}.
}

\noindent{\it Step\,\,1:\,\,An\,\,auxiliary\,\,LDP.}\\
Let $x\geq 0$ be fixed. In this step we show that
$\{S_n(\kappa_n(x))/((1-r^{-1})^{-1}h(x)a_c^{(n)})\}_{n\in\mathbb N}$
obeys an LDP on $\mathbb R$ with speed $v(n):=(1-r^{-1})^{-1}h(x)a_c^{(n)}$ and rate function $I:=H$.
Note that $S_n(\ell_n(x))$, $\ell_n(x):=(\alpha+(1-r^{-1})^{-1}x)a_c^{(n)}$,
is distributed as $\mathrm{Bin}(n-a_n,\pi_n(\ell_n(x)))$. Note also that
by the super-critical condition and the second relation in \eqref{eq:ptcto0} we have
$a_n=v(n)=o(n)$, and by \eqref{eq:speed} and the definition of $h$ easily follows that
$n\pi_n(\ell_n(x))\sim_e v(n)$. Therefore by Lemma \ref{ldp} we have that $\{S_n(\ell_n(x))/((1-r^{-1})^{-1}h(x)a_c^{(n)})\}_{n\in\mathbb N}$
obeys an LDP on $\mathbb R$ with speed $v(n):=(1-r^{-1})^{-1}h(x)a_c^{(n)}$
and rate function $I:=H$. Since the level sets of $H$ are compacts, the claim of this step follows by e.g.
Theorem 4.2.13 p. 130 in \cite{DZ} if we prove that the processes
\[
\{S_n(\kappa_n(x))/((1-r^{-1})^{-1}h(x)a_c^{(n)})\}_{n\in\mathbb N}\quad\text{and}\quad
\{S_n(\ell_n(x))/((1-r^{-1})^{-1}h(x)a_c^{(n)})\}_{n\in\mathbb N},
\]
are exponentially equivalent i.e., for any $\delta\in\R_+$,
\begin{equation}\label{eq:uterm}
\lim_{n\to\infty}\frac{1}{a_c^{(n)}}\log P\left(|S_n(\kappa_n(x))-S_n(\ell_n(x))|>\delta a_c^{(n)}\right)
=-\infty.
\end{equation}
Let $\delta\in\R_+$ be arbitrarily fixed and let $\eta\in\R_+$ be so small that
$r^{-1}((\alpha(x)+\eta)^{r}-(\alpha(x)-\eta)^{r})<\frac{\delta}{1+\eta}$, where
$\alpha(x):=(\alpha+(1-r^{-1})^{-1}x)$. We have
\begin{align}
&\limsup_{n\to\infty}\frac{1}{a_c^{(n)}}\log P\left(|S_n(\kappa_n(x))-S_n(\ell_n(x))|>\delta a_c^{(n)}\right)\nonumber\\
&=\limsup_{n\to\infty}\frac{1}{a_c^{(n)}}\log P\left(\sum_{i\notin\mathcal{A}_n(0)}\ind\{\kappa_n(x)\wedge\ell_n(x)<Y_i^{(n)}
\leq\kappa_n(x)\vee\ell_n(x)\}>\delta a_c^{(n)}\right)\nonumber\\
&\leq\limsup_{n\to\infty}\frac{1}{a_c^{(n)}}\log P\left(\sum_{i\notin\mathcal{A}_n(0)}\ind\{(\alpha(x)-\eta)a_c^{(n)}<Y_i^{(n)}
\leq(\alpha(x)+\eta)a_c^{(n)}\}>\delta a_c^{(n)}\right)\nonumber\\
&=\limsup_{n\to\infty}\frac{1}{a_c^{(n)}}
\log P\left(\mathrm{Bin}(n-a_n,\Pi_n(x,\eta))>\delta a_c^{(n)}\right),\label{eq:limsup}
\end{align}
where
\[
\Pi_n(x,\eta):=\pi_n((\alpha(x)+\eta)a_c^{(n)})-\pi_n((\alpha(x)-\eta)a_c^{(n)}).
\]
For any $n\in\mathbb N$, we clearly have $E[\mathrm{Bin}(n-a_n,\Pi_n(x,\eta))]\leq n\Pi_n(x,\eta)$
and using \eqref{eq:speed} we get
\begin{align}
\limsup_{n\to\infty}\frac{E[\mathrm{Bin}(n-a_n,\Pi_n(x,\eta))]}
{r^{-1}((\alpha(x)+\eta)^{r}-(\alpha(x)-\eta)^{r})a_c^{(n)}}\leq (1-r^{-1})^{r-1}<1.\nonumber
\end{align}
Therefore, by the choice of $\eta$, for all $n$ large enough, we deduce
\[
E[\mathrm{Bin}(n-a_n,\Pi_n(x,\eta)]\leq (1+\eta)r^{-1}((\alpha(x)+\eta)^{r}-(\alpha(x)-\eta)^{r})a_c^{(n)}<\delta a_c^{(n)}.
\]
So, by \eqref{Penrose1.5}, for all $n$ large enough,
\begin{equation}\label{eq:conc}
P\left(\mathrm{Bin}(n-a_n,\Pi_n(x,\eta))>\delta a_c^{(n)}\right)\leq\mathrm{e}^{-(n-a_n)\Pi_n(x,\eta)H\left(\frac{\delta a_c^{(n)}}{(n-a_n)\Pi_n(x,\eta)}\right)}.
\end{equation}
By \eqref{eq:limsup}, \eqref{eq:conc}, \eqref{eq:speed} and $a_n/n\to 0$, we get
\begin{align}
&\limsup_{n\to\infty}\frac{1}{a_c^{(n)}}\log P\left(|S_n(\kappa_n(x))-S_n(\ell_n(x))|>\delta a_c^{(n)}\right)\nonumber\\
&\,\,\,\leq-\frac{1}{r}\left(1-\frac{1}{r}\right)^{r-1}((\alpha(x)+\eta)^{r}-(\alpha(x)-\eta)^{r})
H\left(\frac{\delta}{r^{-1}(1-r^{-1})^{r-1}[(\alpha(x)+\eta)^{r}-(\alpha(x)-\eta)^{r}]}\right).\nonumber
\end{align}
Letting $\eta$ tend to zero we deduce \eqref{eq:uterm} (indeed, $x H(\delta/x)\to +\infty$ as $x\to 0$).
\\
\noindent{\it Step\,\,2:\,\,Proof\,\,of\,\,\eqref{eq:compasB1}.}\\
For technical reasons which will be clear later on, we fix
\[
K>(\alpha+r(r-1)^{-1}x_0)\vee 2\vee (\mathrm{e}^3r(1-r^{-1})^{-(r-1)})^{(r-1.5)^{-1}}
\]
and rewrite the event $B_1^{(n)}$ as the union of the events
\begin{equation*}
B_1^{(n)'}:=\bigcup_{t=a_n}^{\lfloor\kappa_n(x_0)\rfloor}\{S_n(t)+a_n-t\leq 0\}
\end{equation*}
and
\[
B_1^{(n)''}:=\bigcup_{t=\lfloor\kappa_n(x_0)\rfloor}^{\lceil K a_c^{(n)}\rceil}\{S_n(t)+a_n-t\leq 0\}.
\]
We shall show later on
\begin{equation}\label{eq:exacty}
\lim_{n\to\infty}\frac{1}{a_c^{(n)}}\log P(S_n(\kappa_n(x))+a_n-\lfloor\kappa_n(x)\rfloor\leq 0)
=-J(x),\quad x\geq 0.
\end{equation}
Since
\[
P(B_1^{(n)})\geq P(S_n(\kappa_n(x_0))+a_n-\lfloor\kappa_n(x_0)\rfloor\leq 0),\nonumber
\]
by \eqref{eq:exacty} we have
\begin{align}
\liminf_{n\to\infty}\frac{1}{a_c^{(n)}}\log P(B_1^{(n)})&\geq-J(x_0).\label{eq:lower}
\end{align}
Let $t_n\in\{a_n,\ldots,\lfloor\kappa_n(x_0)\rfloor\}$ be such that
\[
\max_{t\in\{a_n,\ldots,\lfloor\kappa_n(x_0)\rfloor\}}P(S_n(t)\leq t-a_n)=P(S_n(t_n)\leq t_n-a_n).
\]
We have
\begin{equation}\label{eq:ineqcorretta}
\limsup_{n\to\infty}\frac{1}{a_c^{(n)}}\log P(S_n(t_n)\leq t_n-a_n)\leq-J(x_0).
\end{equation}
Indeed, reasoning by contradiction, suppose
\begin{equation}\label{eq:reascontra}
\limsup_{n\to\infty}\frac{1}{a_c^{(n)}}\log P(S_n(t_n)\leq t_n-a_n)>-J(x_0).
\end{equation}
Letting $\{t_{n_j}\}_{j\in\mathbb N}$ denote a subsequence of $\{t_n\}_{n\in\mathbb N}$ which realizes this $\limsup$, and setting
\[
x_{n_j}:=\frac{t_{n_j}-\alpha a_c^{(n_j)}}{(1-r^{-1})^{-1}a_c^{(n_j)}}
\]
we have
\[
\frac{a_{n_j}-\alpha a_c^{(n_j)}}{(1-r^{-1})^{-1}a_c^{(n_j)}}\leq x_{n_j}\leq\frac{\lfloor\kappa_{n_j}(x_0)\rfloor-\alpha a_c^{(n_j)}}{(1-r^{-1})^{-1}a_c^{(n_j)}}.
\]
Therefore, by the definition of $\kappa_n(x_0)$,
\[
0\leq\liminf_{j\to\infty}x_{n_j}\leq\limsup_{j\to\infty}x_{n_j}\leq x_0.
\]
So, we may select a subsequence $\{x_{n_{jh}}\}_{h\in\mathbb N}\subseteq\{x_{n_j}\}_{j\in\mathbb N}$
such that $x_{n_{jh}}\to\bar x\in [0,x_0]$, as $h\to\infty$. Consequently, for any $h\in\mathbb N$,
\[
t_{n_{jh}}=(\alpha+(1-r^{-1})^{-1}\bar x+o(1))a_c^{(n_{jh})}=\kappa_{n_{jh}}(\bar x),
\]
and so by \eqref{eq:exacty}
\[
\lim_{h\to\infty}\frac{1}{a_c^{(n_{jh})}}\log P(S_{n_{jh}}(t_{n_{jh}})\leq t_{n_{jh}}-a_{n_{jh}})=-J(\bar x)\leq-J(x_0),
\]
where the latter inequality follows by Lemma \ref{lemma-Hy}$(ii)$. This contradicts \eqref{eq:reascontra} and
proves \eqref{eq:ineqcorretta}, which yields
\begin{align}
\limsup_{n\to\infty}\frac{1}{a_c^{(n)}}\log P(B_1^{(n)'})&\leq\limsup_{n\to\infty}\frac{1}{a_c^{(n)}}\log\left(\sum_{t=a_n}^{\lfloor\kappa_n(x_0)\rfloor}P(S_n(t)\leq t-a_n)\right)\nonumber\\
&\leq\limsup_{n\to\infty}\frac{1}{a_c^{(n)}}\log((\lfloor\kappa_n(x_0)\rfloor-a_n+1)P(S_n(t_n)\leq t_n-a_n))\nonumber\\
&\leq-J(x_0).\label{eq:uB}
\end{align}
Arguing similarly (with obvious modifications), one may check
\begin{align}
\limsup_{n\to\infty}\frac{1}{a_c^{(n)}}\log P(B_1^{(n)''})&\leq-J(x_0).\label{eq:uBbis}
\end{align}
The matching upper bound for \eqref{eq:lower} (and so \eqref{eq:compasB1}) easily follows by the union bound, the principle of the largest term, \eqref{eq:uB} and
\eqref{eq:uBbis}.
To conclude this step, it remains to show \eqref{eq:exacty}. We distinguish two cases: $x>0$ and $x=0$.\\
\noindent{\it Case\,\,1:\,\,$x>0$.} By the super-critical condition
\begin{align}
P(S_n(\kappa_n(x))\leq\lfloor\kappa_n(x)\rfloor-a_n)&=P(S_n(\kappa_n(x))\leq ((1-r^{-1})^{-1}x+o(1))a_c^{(n)})\nonumber\\
&=P\left(\frac{S_n(\kappa_n(x))}{h(x)a_c^{(n)}}\leq\frac{((1-r^{-1})^{-1}x+o(1))}{h(x)}\right)\nonumber\\
&=P\left(\frac{S_n(\kappa_n(x))}{(1-r^{-1})^{-1}h(x)a_c^{(n)}}\leq\frac{(x+o(1))}{h(x)}\right).\label{eq:AAPATT}
\end{align}
So for $\varepsilon$ arbitrarily chosen in $(0,x\wedge h(x))\equiv (0,x)$ (see Lemma \ref{lemma-Hy}$(i)$)
and $n$ large enough
\begin{align}
\log P(S_n(\kappa_n(x))\leq\lfloor\kappa_n(x)\rfloor-a_n)\geq
\log P\left(\frac{S_n(\kappa_n(x))}{(1-r^{-1})^{-1}h(x)a_c^{(n)}}<\frac{\varepsilon}{h(x)}\right).\label{eq:19}
\end{align}
By the LDP in Step 1 we have
\begin{align}
\liminf_{n\to\infty}\frac{1}{(1-r^{-1})^{-1}h(x)a_c^{(n)}}\log P\left(\frac{S_n(\kappa_n(x))}{(1-r^{-1})^{-1}h(x)a_c^{(n)}}<\frac{\varepsilon}{h(x)}\right)
&\geq-\inf_{y\in\left(-\infty,\frac{\varepsilon}{h(x)}\right)}H(y)\nonumber\\
&=-H\left(\frac{\varepsilon}{h(x)}\right)\label{eq:20},
\end{align}
where in \eqref{eq:20} we used that $H\equiv +\infty$ on $\R_-$ and that $H$ is continuously decreasing on $[0,1)$. By \eqref{eq:19} and \eqref{eq:20}, we deduce
\begin{align}
\liminf_{n\to\infty}\frac{1}{a_c^{(n)}}\log P(S_n(\kappa_n(x))\leq\lfloor\kappa_n(x)\rfloor-a_n)&
\geq-(1-r^{-1})^{-1}h(x)H\left(\frac{\varepsilon}{h(x)}\right)
\nonumber
\end{align}
and taking the supremum over $\varepsilon\in (0,x)$ by the properties of $H$ and Lemma \ref{lemma-Hy}$(i)$ we get the lower bound
\begin{align}
\liminf_{n\to\infty}\frac{1}{a_c^{(n)}}\log P(S_n(\kappa_n(x))\leq\lfloor\kappa_n(x)\rfloor-a_n)&
\geq-J(x).
\label{eq:lowerHK}
\end{align}
Now we prove the matching upper bound (we remark that the proof of the matching upper bound we are going to give still holds for $x=0$).
For $\varepsilon$ arbitrarily chosen in $(x,\infty)$ and $n$ enough,
by \eqref{eq:AAPATT} and the LDP in Step 1 we have
\begin{align}
\limsup_{n\to\infty}\frac{1}{a_c^{(n)}}\log P(S_n(\kappa_n(x))\leq\lfloor\kappa_n(x)\rfloor-a_n)&
\leq\limsup_{n\to\infty}\frac{1}{a_c^{(n)}}\log P\left(\frac{S_n(\kappa_n(x))}{(1-r^{-1})^{-1}h(x)a_c^{(n)}}\leq
\frac{\varepsilon}{h(x)}\right)\nonumber\\
&\leq-(1-r^{-1})^{-1}h(x)\inf_{y\in \left(-\infty,
\frac{\varepsilon}{h(x)}
\right]}H(y).\nonumber
\end{align}
Taking the infimum over $\varepsilon$ we then have
\begin{align}
\limsup_{n\to\infty}\frac{1}{a_c^{(n)}}\log P(S_n(\kappa_n(x))\leq\lfloor\kappa_n(x)\rfloor-a_n)
&\leq-(1-r^{-1})^{-1}h(x)\sup_{\varepsilon>x}
\inf_{y\in \left(-\infty,
\frac{\varepsilon}{h(x)}
\right]}H(y)\nonumber\\
&=-(1-r^{-1})^{-1}h(x)
\inf_{y\in \left(-\infty,
\frac{x}{h(x)}
\right]}H(y)\nonumber\\
&=-J(x),\label{eq:up}
\end{align}
where the latter equality is a consequence of Lemma \ref{lemma-Hy}$(i)$ and the fact that $H\equiv +\infty$ on $\R_-$ and $H$ is continuously decreasing on $[0,1)$.
Relation \eqref{eq:exacty} follows by \eqref{eq:lowerHK} and \eqref{eq:up}.
\\
\noindent{\it Case\,\,2:\,\,$x=0$.}
We have
\begin{align}
P(S_n(\kappa_n(0))\leq\lfloor\kappa_n(0)\rfloor-a_n)&=P\left(\frac{S_n(\kappa_n(0))}{a_c^{(n)}}\leq o(1)\right)\nonumber\\
&\geq P\left(S_n(\kappa_n(0))=0\right)=(1-\pi_n(\kappa_n(0)))^{n-a_n}.\label{eq:primo}
\end{align}
Therefore
\begin{align}
\liminf_{n\to\infty}\frac{1}{a_c^{(n)}}\log P(S_n(\kappa_n(0))\leq\lfloor\kappa_n(0)\rfloor-a_n)&\geq
\liminf_{n\to\infty}\frac{n-a_n}{a_c^{(n)}}\log(1-\pi_n(\kappa_n(0)))\nonumber\\
&=\liminf_{n\to\infty}\left(-\frac{n}{a_c^{(n)}}\pi_n(\kappa_n(0))+\frac{n}{a_c^{(n)}}o(\pi_n(\kappa_n(0)))\right)\nonumber\\
&=-r^{-1}(1-r^{-1})^{r-1}\alpha^r=-J(0).
\end{align}
where we used \eqref{eq:speed} (which yields $\pi_n(\kappa_n(0))\sim_e r^{-1}(1-r^{-1})^{r-1}\alpha^r a_c^{(n)}/n$).
The proof of the matching upper bound has been already done (see the Case $1$).
\\
\noindent{\it Step\,\,3:\,\,Proof\,\,of\,\,\eqref{eq:compasb2}.}\\
For $n\geq 2$
define
\[
J_n:=\min\{j\geq\lceil K\rceil:\,\,p_n \theta_j^{(n)}\geq 1\},\quad\text{where}\quad\theta_j^{(n)}:=K^{j/\lceil K\rceil}a_c^{(n)}.
\]
By construction we have
\begin{equation}\label{eq:construction}
\theta_{J_n}^{(n)}\geq\lfloor p_n^{-1}\rfloor,\quad
p_n\theta_{J_n}^{(n)}
<
K^\frac{1}{K}\leq 2,
\end{equation}
where
the second inequality is a consequence of the relations $\theta_j^{(n)}=\theta_{j-1}^{(n)}K^{1/\lceil K\rceil}$ and $p_n\theta_{J_{n-1}}^{(n)}<1$. By \eqref{eq:dis2} and
\eqref{eq:construction}
\begin{align}
P(B_2^{(n)})
&\leq\sum_{j=\lceil K\rceil}^{J_n-1}P\left(\bigcup_{t\in [\theta_j^{(n)},\theta_{j+1}^{(n)}]\cap\mathbb N}\{S_n'(t)\leq t\}\right)
\leq\sum_{j=\lceil K\rceil}^{J_n-1}P\left(S_n'(\theta_j^{(n)})\leq\theta_{j+1}^{(n)}\right),\label{eq:psprimo}
\end{align}
where for the latter inequality we used that $S_n'(t)$ is non-decreasing with respect to $t$ and that
$\theta_j^{(n)}$ is non-decreasing with respect to $j$ (this latter monotonicity is guaranteed by the fact that $K>1$).
For all $n$ large enough and $j\in\{\lceil K\rceil,\ldots,J_n-1\}$, by the usual Poisson approximation for the binomial distribution
we have
\begin{align}
\pi_n(\theta_j^{(n)})&\geq
\frac{(p_n\lfloor\theta_j^{(n)}\rfloor)^{r}}{r!}\mathrm{e}^{-p_n\lfloor\theta_j^{(n)}\rfloor}(1+o(1))\nonumber\\
&>(1+o(1))\mathrm{e}^{-2}\frac{(p_n\lfloor\theta_j^{(n)}\rfloor)^{r}}{r!}\label{eq:ptczero}\\
&=(1+o(1))\mathrm{e}^{-2}(1-r^{-1})^{r-1}\frac{\lfloor\theta_j^{(n)}\rfloor}{nr}\left(\frac{\lfloor\theta_j^{(n)}\rfloor}{a_c^{(n)}}\right)^{r-1}\label{eq:deftc}\\
&\geq\mathrm{e}^{-3}(1-r^{-1})^{r-1}\frac{(K^{j/\lceil K\rceil})^{r-1}}{nr}\theta_j^{(n)},
\label{eq:pitj}
\end{align}
where in \eqref{eq:ptczero} we used the second relation in \eqref{eq:construction} and in \eqref{eq:deftc} we used the definition of $a_c^{(n)}$.
Therefore,
for $n$ large enough, we deduce
\begin{align}
E[S_n'(\theta_j^{(n)})]=n\pi_n(\theta_j^{(n)})
&\geq
\mathrm{e}^{-3}(1-r^{-1})^{r-1}\frac{K^{r-1}}{r}\theta_j^{(n)}>K^{1/\lceil K\rceil}\theta_j^{(n)}=\theta_{j+1}^{(n)}.\label{eq:meansprimo}
\end{align}
So, by \eqref{Penrose1.6}, for all $n$ large enough,
\[
P(S_n'(\theta_j^{(n)})\leq\theta_{j+1}^{(n)})\leq\mathrm{e}^{-n\pi_n(\theta_j^{(n)})H(\theta_{j+1}^{(n)}/(n\pi_n(\theta_j^{(n)})))},\quad\text{$j=\lceil K\rceil,\ldots,J_n-1$.}
\]
By \eqref{eq:meansprimo}, for $n$ large enough and $j=\lceil K\rceil,\ldots,J_n-1$,
\[
1>x_K:=\frac{r K^{1/\lceil K\rceil}}{\mathrm{e}^{-3}(1-r^{-1})^{r-1}K^{r-1}}
=\frac{\theta_{j+1}^{(n)}}{\mathrm{e}^{-3}(1-r^{-1})^{r-1}\frac{K^{r-1}}{r}\theta_j^{(n)}}\geq\frac{\theta_{j+1}^{(n)}}{n\pi_n(\theta_j^{(n)})}>0.
\]
Therefore, using that $H$ is decreasing on $[0,1)$ and \eqref{eq:pitj},
\begin{align}
\pi_n(\theta_j^{(n)})H\left(\frac{\theta_{j+1}^{(n)}}{n\pi_n(\theta_j^{(n)})}\right)&\geq\pi_n(\theta_j^{(n)})H(x_K)\geq
\mathrm{e}^{-3}(1-r^{-1})^{r-1}\frac{(K^{r/\lceil K\rceil})^j}{nr}
H(x_K)a_c^{(n)},\nonumber
\end{align}
and so, for $n$ large enough and $j=\lceil K\rceil,\ldots,J_n-1$,
\begin{align}
P(S_n'(\theta_j^{(n)})\leq\theta_{j+1}^{(n)})
&\leq\mathrm{e}^{-\mathrm{e}^{-3}r^{-1}(1-r^{-1})^{r-1}K^{r}H(x_K)a_c^{(n)}}
\mathrm{e}^{-\left(K^{\frac{(j-\lceil K\rceil)r}{\lceil K\rceil}}-1\right)\mathrm{e}^{-3}r^{-1}(1-r^{-1})^{r-1}K^{r}H(x_K)a_c^{(n)}}\nonumber\\
&\leq\mathrm{e}^{-\mathrm{e}^{-3}r^{-1}(1-r^{-1})^{r-1}K^{r}H(x_K)a_c^{(n)}}
\mathrm{e}^{-(\log K)(j-\lceil K\rceil)(\lceil K\rceil)^{-1}\mathrm{e}^{-3}(1-r^{-1})^{r-1}K^{r}H(x_K)a_c^{(n)}}
\label{eq:disugetater}
\end{align}
where
we used the relation $K^{\frac{(j-\lceil K\rceil)r}{\lceil K\rceil}}-1\geq\frac{(j-\lceil K\rceil)r}{\lceil K\rceil}\log K$.
By \eqref{eq:psprimo} and \eqref{eq:disugetater}, we have
\begin{align}
P(B_2^{(n)})&\le
\mathrm{e}^{-\mathrm{e}^{-3}r^{-1}(1-r^{-1})^{r-1}K^{r}H(x_K)a_c^{(n)}}
\frac{1}{1-\mathrm{e}^{-\mathrm{e}^{-3}(1-r^{-1})^{r-1}K^{r}(\lceil K\rceil)^{-1}H(x_K)(\log K) a_c^{(n)}}}.\nonumber
\end{align}
Relation \eqref{eq:compasb2} follows by this inequality setting
$C_K:=\frac{\mathrm{e}^{-3}(1-r^{-1})^{r-1}K^{r}H(x_{\lceil K\rceil})}{r}>0$
and choosing $K$ so large that $C_K>J(x_0)$.
\\
\noindent{\it Step\,\,4:\,\,Proof\,\,of\,\,\eqref{eq:compasb3}.}\\
For $n$ large enough, we have
\begin{align}
\pi_n(p_n^{-1})
&=P(\mathrm{Bin}(\lfloor p_n^{-1}\rfloor,p_n)\geq r)=P(\mathrm{Po}(\lfloor p_n^{-1}\rfloor p_n)\geq r)+O(p_n)\geq 2c,\label{eq:b3c}
\end{align}
for some small $c\in (0,1)$, see e.g. the proof of Lemma 8.2 Case 3 p. 26 in \cite{JLTV}.
By \eqref{eq:dis2}, for $n$ sufficiently large, we have
\begin{align}
B_3^{(n)}&\subseteq
\bigcup_{t=\lfloor p_n^{-1}\rfloor}^{\lfloor cn\rfloor}\{S_n'(t)\leq t\}\cup\bigcup_{t=\lfloor cn\rfloor}^{\lfloor n-p_n^{-1}\rfloor\wedge\lfloor n-\varepsilon f_1(n)\rfloor}\{S_n'(t)\leq t\}
\cup\bigcup_{t=\lfloor n-p_n^{-1}\rfloor\wedge\lfloor n-\varepsilon f_1(n)\rfloor}^{\lfloor n-\varepsilon f_1(n)\rfloor}\{S_n'(t)\leq t\},
\nonumber
\end{align}
with the convention that the latter union of events is empty if $\lfloor n-p_n^{-1}\rfloor>\lfloor n-\varepsilon f_1(n)\rfloor$ for all $n$ large enough.
From now on, we suppose $f_1(n)=o(p_n^{-1})$. The case $p_n f_1(n)\to\ell\in (0,\infty]$ may be treated with obvious modifications. If $f_1(n)=o(p_n^{-1})$, then
\begin{align}
B_3^{(n)}&\subseteq
\bigcup_{t=\lfloor p_n^{-1}\rfloor}^{\lfloor cn\rfloor}\{S_n'(t)\leq t\}\cup\bigcup_{t=\lfloor cn\rfloor}^{\lfloor n-p_n^{-1}\rfloor}\{S_n'(t)\leq t\}
\cup\bigcup_{t=\lfloor n-p_n^{-1}\rfloor}^{\lfloor n-\varepsilon f_1(n)\rfloor}\{S_n'(t)\leq t\}
\label{eq:hypfop}\\
&\subseteq\{S_n'(p_n^{-1})\leq \lfloor cn \rfloor\}\cup\{S_n'(cn)\leq\lfloor n-p_n^{-1}\rfloor\}
\cup\{S_n'(n-p_n^{-1})\leq\lfloor n-\varepsilon f_1(n)\rfloor\},\label{eq:upperb3n}
\end{align}
where for the latter inclusion we used that the events $\{S_n'(u)\leq v\}$ are non-increasing in $u$ and non-decreasing in $v$.
By \eqref{eq:b3c}, for $n$ large enough,
\[
cn\leq E[S_n'(p_n^{-1})]=n\pi_n(p_n^{-1}).
\]
So by \eqref{Penrose1.6}, for $n$ large enough,
\begin{align}
P(S_n'(p_n^{-1})\leq\lfloor cn\rfloor)&\leq\mathrm{e}^{-n\pi_n(p_n^{-1})H(\lfloor cn\rfloor /(n\pi_n(p_n^{-1})))}\nonumber\\
&\leq\mathrm{e}^{-n\pi_n(p_n^{-1})H(c/\pi_n(p_n^{-1}))}\label{eq:H01:1}\\
&\leq\mathrm{e}^{-cn H(c/\pi_n(p_n^{-1}))}\nonumber\\
&\leq\mathrm{e}^{-cH(1/2) n},\label{eq:H01:2}
\end{align}
where the inequalities \eqref{eq:H01:1} and \eqref{eq:H01:2} follow recalling that $H$ decreases on $[0,1]$ and using \eqref{eq:b3c}. Therefore, by the second limit
in \eqref{eq:ptcto0}, we have
\begin{align}
\limsup_{n\to\infty}\frac{1}{a_c^{(n)}}\log P(S_n'(p_n^{-1})\leq\lfloor cn\rfloor)&\leq-cH(1/2)\lim_{n\to\infty}\frac{n}{a_c^{(n)}}=-\infty.\label{eq:dis3}
\end{align}
By \eqref{eq:hyp} we deduce $n p_n\to\infty$. So, using again \eqref{Penrose1.6}, for $n$ large enough,
\begin{align}
\pi_n(cn)&=P(\mathrm{Bin}(\lfloor cn\rfloor,p_n)\geq r)\geq 1-\mathrm{e}^{-\lfloor cn\rfloor p_n H(r/(\lfloor cn\rfloor p_n))};\label{eq:lim0}
\end{align}
moreover,
\begin{equation}\label{eq:asrH}
\lfloor cn\rfloor p_n H(r/(\lfloor cn\rfloor p_n))\sim_e\lfloor cn\rfloor p_n.
\end{equation}
By \eqref{eq:lim0} and \eqref{eq:asrH} we have
\begin{equation}\label{eq:limite}
\lim_{n\to\infty}\frac{1}{np_n(1-\pi_n(cn))}=+\infty.
\end{equation}
Therefore, applying \eqref{Penrose1.7}, we deduce
\begin{align}
P(S_n'(cn)\leq \lfloor n-p_n^{-1}\rfloor)&=P(n-S_n'(cn)\geq n-\lfloor n-p_n^{-1}\rfloor)\nonumber\\
&=P(\mathrm{Bin}(n,1-\pi_n(cn))\geq n-\lfloor n-p_n^{-1}\rfloor)\nonumber\\
&=P(\mathrm{Bin}(n,1-\pi_n(cn))\geq\lceil p_n^{-1}\rceil)\label{eq:nmeno}\\
&\leq\mathrm{e}^{-\frac{\lceil p_n^{-1}\rceil}{2}\log\left(\frac{\lceil p_n^{-1}\rceil}{n(1-\pi_n(cn))}\right)}.\nonumber
\end{align}
Here in \eqref{eq:nmeno} we used that, for any $n\in\mathbb N$ and $y\in\mathbb R$,
\begin{equation}\label{eq:intpart}
n-\lfloor n-y\rfloor=n+\lceil y-n\rceil=n+\lceil y\rceil-n=\lceil y\rceil.
\end{equation}
By this latter inequality, \eqref{eq:limite} and the third limit in \eqref{eq:ptcto0}, we have
\begin{align}
\limsup_{n\to\infty}\frac{1}{a_c^{(n)}}\log P(S_n'(cn)\leq \lfloor n-p_n^{-1}\rfloor)=
-\infty.\label{eq:-infty}
\end{align}
Arguing as for \eqref{eq:nmeno}, we have
\begin{align}
P(S_n'(n-p_n^{-1})\leq \lfloor n-\varepsilon f_1(n)\rfloor)
&=P(\mathrm{Bin}(n,1-\pi_n(n-p_n^{-1}))\geq\lceil\varepsilon f_1(n)\rceil).\label{eq:nmenobis}
\end{align}
Let $\delta\in\mathbb R_+$ be so small that $\mathrm{e}^{-1}+\delta<1$. By \eqref{eq:svilpi}
\begin{align}
\lim_{n\to\infty}\frac{f_1(n)}{n(1-\pi_n(n-p_n^{-1}))}
&\geq\lim_{n\to\infty}\frac{f_1(n)}{n(n p_n)^{r-1}[(1-p_n)^{p_n^{-1}}]^{n p_n-1}}\nonumber\\
&\geq\mathrm{e}^{-1}\lim_{n\to\infty}\frac{a_c^{(n)}g(n)}{n (n p_n)^{r}(\mathrm{e}^{-1}+\delta)^{n p_n}}=+\infty,\nonumber
\end{align}
where the latter relation follows by the definition of $f_1$ and noticing that by \eqref{eq:tc}
\begin{equation*}
\frac{a_c^{(n)}g(n)}{n (n p_n)^{r}(\mathrm{e}^{-1}+\delta)^{n p_n}}\sim g(n)(np_n)^{-r^2/(r-1)}(\mathrm{e}^{-1}+\delta)^{-np_n}\to\infty.
\end{equation*}
Therefore, by applying again \eqref{Penrose1.7}
\begin{align}
P(S_n'(n-p_n^{-1})\leq \lfloor n-\varepsilon f_1(n)\rfloor)&\leq\exp\left(-\frac{\lceil\varepsilon f_1(n)\rceil}{2}\log\left(\frac{\lceil\varepsilon f_1(n)\rceil}{n(1-\pi_n( n-p_n^{-1}))}\right)\right).\nonumber
\end{align}
Consequently, by the definition of $f_1$ and \eqref{eq:svilpi}
\begin{align}
&\limsup_{n\to\infty}\frac{1}{a_c^{(n)}}\log P(S_n'(n-p_n^{-1})\leq \lfloor n-\varepsilon f_1(n)\rfloor)\leq-\frac{\varepsilon}{2}
\limsup_{n\to\infty}\frac{f_1(n)}{a_c^{(n)}}\log\left(\frac{f_1(n)}{n(n p_n)^{r-1}(1-p_n)^{n-p_n^{-1}}}
\right)\nonumber\\
&\,\,\,\,\,\,
\,\,\,\,\,\,
\,\,\,\,\,\,
\,\,\,\,\,\,
\leq-\frac{\varepsilon}{2}
\limsup_{n\to\infty}\frac{g(n)}{np_n}\log\left(\frac{a_c^{(n)}g(n)}{n(n p_n)^{r}(1-p_n)^{n-p_n^{-1}}}\right)\nonumber\\
&\,\,\,\,\,\,
\,\,\,\,\,\,
\,\,\,\,\,\,
\,\,\,\,\,\,
=-\frac{\varepsilon}{2}
\limsup_{n\to\infty}\frac{g(n)}{np_n}\log(g(n)(np_n)^{-r^2/(r-1)}(1-p_n)^{-n+p_n^{-1}})\nonumber\\
&\,\,\,\,\,\,
\,\,\,\,\,\,
\,\,\,\,\,\,
\,\,\,\,\,\,
\leq-\frac{\varepsilon}{2}
\lim_{n\to\infty}g(n)\lim_{n\to\infty}\frac{1}{np_n}\left(np_n-\frac{r^2}{r-1}\log(np_n)\right)=-\infty.\nonumber
\end{align}
The claim follows by this latter relation,
\eqref{eq:upperb3n}, \eqref{eq:dis3}, \eqref{eq:-infty} and the principle of the largest term.\\

\subsubsection{Proof of Theorem \ref{ldp1}}
We divide the proof in two steps. In the first step we prove the large deviation lower bound and
in the second step we prove the large deviation upper bound.\\
\noindent{\it Step\,\,1:\,\,large\,\,deviation\,\,lower\,\,bound.}\\
Let $O\subseteq\overline{\mathbb R}$ be an open set. If $0\in O$, then since $O$ is open there exists $\delta>0$ such that $(-\delta,\delta)\subset O$.
For a fixed $0<\varepsilon<\delta\wedge\ell_1^{-1}$, by \eqref{eq:ratiopos} we have
\begin{align}
P\left(\frac{n-A_n^*}{f_1(n)}\in O\right)\geq P\left(\frac{n-A_n^*}{f_1(n)}\in [0,\delta)\right)
\geq P\left(\frac{n-A_n^*}{f_1(n)}\in [0,\varepsilon]\right),\quad n\in\mathbb N.\label{eq:disin}
\end{align}
By Lemma \ref{lemma-Hy} we have $J(x_0)>0$. Let $0<\eta<J(x_0)$ be arbitrarily fixed. By Proposition \ref{ldp:supcrit} we have that there exists $n_\eta$ such that, for any $n>n_\eta$,
\begin{equation}\label{eq:exp1}
1-\mathrm{e}^{-(J(x_0)-\eta)a_c^{(n)}}<P\left(\frac{n-A_n^*}{f_1(n)}\leq \varepsilon\right)<1-\mathrm{e}^{-(J(x_0)+\eta)a_c^{(n)}}.
\end{equation}
By \eqref{eq:disin} and \eqref{eq:exp1}, we easily have
\[
\liminf_{n\to\infty}\frac{1}{a_c^{(n)}}\log P\left(\frac{n-A_n^*}{f_1(n)}\in O\right)\geq 0=-\inf_{x\in O} I_1(x),
\]
and the large deviation lower bound for the case $0\in O$ is proved. If $0\notin O$, then the claim is obvious if in addition $\ell_1^{-1}\notin O$. Otherwise, we distinguish
two further cases: $\ell_1^{-1}=+\infty\in O$ or $\ell_1^{-1}\in O\cap (0,\infty)$. If $\ell_1^{-1}=+\infty\in O$, then there exists $\varepsilon>0$ such that $O\supseteq (\varepsilon,+\infty]$.
So by Proposition \ref{ldp:supcrit}
\begin{align}
\liminf_{n\to\infty}\frac{1}{a_c^{(n)}}\log P\left(\frac{n-A_n^*}{f_1(n)}\in O\right)&\geq
\liminf_{n\to\infty}\frac{1}{a_c^{(n)}}\log P\left(\frac{n-A_n^*}{f_1(n)}>\varepsilon\right)\nonumber\\
&=-J(x_0)=-\inf_{x\in O}I_1(x).\nonumber
\end{align}
If $\ell_1^{-1}\in O\cap (0,\infty)$, then since $O\cap (0,\infty)$ is open, there exists $\delta\in (0,\ell_1^{-1})$ such that $(\ell_1^{-1}-\delta,\ell_1^{-1}+\delta]\subset O\cap (0,\infty)$. Therefore,
for all $n$ large enough,
\begin{align}
P\left(\frac{n-A_n^*}{f_1(n)}\in O\right)&\geq P\left(\frac{n-A_n^*}{f_1(n)}\in (\ell_1^{-1}-\delta,\ell_1^{-1}+\delta]\right)\nonumber\\
&=P\left(\frac{n-A_n^*}{f_1(n)}\leq\ell_1^{-1}+\delta\right)-P\left(\frac{n-A_n^*}{f_1(n)}\leq\ell_1^{-1}-\delta\right)\nonumber\\
&=1-P\left(\frac{n-A_n^*}{f_1(n)}\leq\ell_1^{-1}-\delta\right)\nonumber\\
&=P\left(\frac{n-A_n^*}{f_1(n)}>\ell_1^{-1}-\delta\right),\nonumber
\end{align}
where we used that, for $n$ large enough,
\[
\frac{n-A_n^*}{f_1(n)}\leq\frac{n}{f_1(n)}<\ell_1^{-1}+\delta.
\]
The large deviation lower bound easily follows by Proposition \ref{ldp:supcrit}.\\
\noindent{\it Step\,\,2:\,\,large\,\,deviation\,\,upper\,\,bound.}\\
Let $C\subseteq\overline{\mathbb R}$ be a closed set.
If $0\in C$, then the large deviation upper bound is trivial. If $0\notin C$, we start noticing that by \eqref{eq:ratiopos}
and $n/f_1(n)\to\ell_1^{-1}$, for any $\delta\in\mathbb R_+$ and all $n$ large enough,
\begin{equation}\label{eq:upperclosed}
P\left(\frac{n-A_n^*}{f_1(n)}\in C\right)=P\left(\frac{n-A_n^*}{f_1(n)}\in C\cap (0,\ell_1^{-1}+\delta)\right).
\end{equation}
Then we distinguish two cases: $\ell_1^{-1}\in C$ and $\ell_1^{-1}\notin C$.
If $\ell_1^{-1}\in C$, then there exists $\varepsilon\in (0,\ell_1^{-1})$ such that $C\cap (0,\ell_1^{-1}+\delta)\subset (\varepsilon,+\infty]$
and the large deviation upper bound easily follows by \eqref{eq:upperclosed} and Proposition \ref{ldp:supcrit}. If $\ell_1^{-1}\notin C$, then
\[
P\left(\frac{n-A_n^*}{f_1(n)}\in C\right)=P\left(\frac{n-A_n^*}{f_1(n)}\in C\cap (0,\ell_1^{-1})\right)+P\left(\frac{n-A_n^*}{f_1(n)}\in C\cap (\ell_1^{-1},\ell_1^{-1}+\delta)\right).
\]
Since
\[
\lim_{n\to\infty}\frac{1}{a_c^{(n)}}\log P\left(\frac{n-A_n^*}{f_1(n)}\in C\cap (\ell_1^{-1},\ell_1^{-1}+\delta)\right)=-\infty,
\]
recalling that $0,\ell_1^{-1}\notin C$, the large deviation upper bound easily follows by the principle of the largest term
(see e.g. Lemma 1.2.15 p. 7 in \cite{DZ}), if we prove
\[
\lim_{n\to\infty}\frac{1}{a_c^{(n)}}\log P\left(\frac{n-A_n^*}{f_1(n)}\in C\cap (0,\ell_1^{-1})\right)=-\infty.
\]
To this aim, we start noticing that there exist $\kappa_1,\kappa_2\in (0,\ell_1^{-1})$ such that $C\cap (0,\ell_1^{-1})\subset (\kappa_1,\kappa_2]$, and so
\begin{align}
P\left(\frac{n-A_n^*}{f_1(n)}\in C\cap (0,\ell_1^{-1})\right)&\leq P\left(\kappa_1<\frac{n-A_n^*}{f_1(n)}\leq\kappa_2\right)\nonumber\\
&=P\left(\frac{n-A_n^*}{f_1(n)}>\kappa_1\right)-P\left(\frac{n-A_n^*}{f_1(n)}>\kappa_2\right).\label{eq:kappa12}
\end{align}
By this relation and \eqref{eq:pdopo}, we have
\begin{align}
P\left(\frac{n-A_n^*}{f_1(n)}\in C\cap (0,\ell_1^{-1})\right)&\leq
P(B^{(n)}),\nonumber
\end{align}
where
\[
B^{(n)}:=\bigcup_{t=\lfloor n-\kappa_2 f_1(n)\rfloor}^{\lfloor n-\kappa_1 f_1(n)\rfloor}\{S_n(t)+a_n-t\leq 0\},
\]
and so it suffices to prove
\begin{equation}\label{eq:fineprova}
\lim_{n\to\infty}\frac{1}{a_c^{(n)}}\log P(B^{(n)})=-\infty.
\end{equation}
To this aim, we note that
\[
\frac{p_n^{-1}}{n-\kappa_2 f_1(n)}\to 0,
\]
so for $n$ large enough
\begin{equation*}
\bigcup_{t=\lfloor p_n^{-1}\rfloor}^{\lfloor n-\kappa_1 f_1(n)\rfloor}\{S_n(t)+a_n-t\leq 0\}\supset B^{(n)}
\end{equation*}
and consequently \eqref{eq:fineprova} follows by \eqref{eq:compasb3}.

\subsubsection{Proofs of lemmas}

\noindent{\it Proof\,\,of\,\,Lemma\,\,\ref{ldp}.}\\
We shall apply the G\"artner-Ellis Theorem (see e.g. Theorem 2.3.6 p. 44 in \cite{DZ}). Denote by $\mathrm{Be}(\pi(n))$
a Bernoulli distributed random variable with mean $\pi(n)$. For any $\theta\in\mathbb R$, we have
\begin{align}
\Lambda_{n}(\theta):=\log E\left[\mathrm{e}^{\theta\mathrm{Bin}(n-a_n,\pi(n))}\right]
&=(n-a_n)\log E\left[\mathrm{e}^{\theta\mathrm{Be}(\pi(n))}\right]\nonumber\\
&=(n-a_n)\log(1+\pi(n)(\mathrm{e}^\theta-1)).\nonumber
\end{align}
By this relation and the assumptions of the lemma
we deduce
\begin{align}
\lim_{n\to\infty}\frac{\Lambda_n(\theta)}{v(n)}&
=\lim_{n\to\infty}\frac{\log(1+\pi(n)(\mathrm{e}^\theta-1))}{\pi(n)}\lim_{n\to\infty}\frac{n-a_n}{n}\nonumber\\
&=\lim_{n\to\infty}\frac{\mathrm{e}^\theta-1}{1+\pi(n)(\mathrm{e}^\theta-1)}\label{eq:hop}\\
&=\mathrm{e}^{\theta}-1,\nonumber
\end{align}
where \eqref{eq:hop} follows by l'Hopital's rule, which is applicable since $\pi(n)\to 0$.
A straightforward computation shows that the convex conjugate of the function $\theta\mapsto\mathrm{e}^\theta-1$ is $H$.
By the G\"artner-Ellis Theorem we have that $\Biggl\{\frac{\mathrm{Bin}(n-a_n,\pi(n))}{v(n)}\Biggr\}_{n\in\mathbb N}$ satisfies the large deviation
upper bound over the closed sets with speed $v(n)$ and rate function $I:=H$ and, for any open set $O\subseteq\R$, it is satisfied the lower bound
\begin{equation}\label{eq:exposed}
\liminf_{n\to\infty}\frac{1}{v(n)}\log P\left(\frac{\mathrm{Bin}(n-a_n,\pi(n))}{v(n)}\in O\right)\geq-\inf_{y\in O\cap\mathcal F}H(y),
\end{equation}
where $\mathcal F$ is the set of exposed points of $H$ whose exposing hyperplane belongs to $\R$ (we refer to \cite{DZ} for these notions). By Lemma 2.3.9$(b)$ p. 46 in \cite{DZ}
we have $\mathcal F\supseteq\R_+$. Setting $\mathcal F_-:=\mathcal F\cap (-\infty,0]$, by the properties of $H$ we
have
\[
\inf_{y\in O\cap\mathcal F}H(y)=\min\{\inf_{y\in O\cap\mathcal{F}_-}H(y),\inf_{y\in O\cap\R_+}H(y)\}=\inf_{y\in O\cap\R_+}H(y)=\inf_{y\in O}H(y).
\]
The large deviation lower bound over the open sets then follows by this relation and \eqref{eq:exposed}.\\

\noindent{\it Proof\,\,of\,\,Lemma\,\,\ref{lemma-Hy}.}\\
A simple computation shows that the function $x\mapsto\frac{x}{h(x)}$ is strictly increasing on $[0,\alpha/r]$
and strictly decreasing on $(\alpha/r,\infty)$. In particular, $\alpha/r$ is the unique point of maximum of $x\mapsto\frac{x}{h(x)}$ on $[0,\infty)$.
The claim $(i)$ then follows noticing that $\frac{\alpha/r}{h(\alpha/r)}=\frac{1}{\alpha^{r-1}}<1$, whenever $\alpha>1$. We now show $(ii)$.
Since $H$ is strictly decreasing on $[0,1)$ and $x\mapsto\frac{x}{h(x)}$ is strictly decreasing on $(\alpha/r,\infty)$,
we have that $J$ is strictly increasing on $(\alpha/r,\infty)$. For ease of notation, set $\widetilde{h}(x):=x/h(x)$. A simple computation shows that the second derivative
of $(1-r^{-1})J$ is equal to
\[
h''(x)\left(1-\widetilde{h}(x)\right)+\frac{(h(x)\widetilde{h}^{\prime}(x))^2}{x},
\]
which is strictly positive on $(0,\alpha/r]$ and therefore $J$ is strictly convex on $(0,\alpha/r]$.
Furthermore, we note that $\lim_{x\to 0}J'(x)=-\infty$. Collecting all these properties of $J$
we finally deduce $(ii)$.
\\
\noindent$\square$

\subsection{Proofs of Proposition \ref{ldp:supcritbis} and of Theorem \ref{ldp2}}

\subsubsection{Proof of Proposition \ref{ldp:supcritbis}}

\noindent{\it Proof\,\,of\,\,\eqref{eq:lim1}.}\\
Arguing as for \eqref{eq:pdopo}, for any $\varepsilon>\ell_2^{-1}$,
\begin{align}
P\left(\frac{n-A_n^*}{f_2(n)}>\varepsilon\right)&=
P\left(\bigcup_{t=a_n}^{\lfloor n-\varepsilon f_2(n)\rfloor}\{S_n(t)+a_n-t\leq 0\}\right).\nonumber
\end{align}
Let $f_1$ be defined as in the statement of Theorem \ref{ldp1}. Since $b_c^{(n)}\to\infty$, by the third relation in
\eqref{eq:ptcto0} we have $a_c^{(n)}/(n p_n)\to\infty$, and so, for $n$ large enough,
$f_1(n)=(g(n) a_c^{(n)})/(n p_n)$. For later purposes, we choose $g$ in such a way that $\ell_1\in [0,1)$ and $(g(n)a_c^{(n)})/n\to 0$
(so that $f_1(n)=o(p_n^{-1})$). Note that, by the latter relation in \eqref{eq:ptcto0} and $f_2(n)\sim\ell_2 b_c^{(n)}$, for any fixed positive constant $K$,
$f_1(n)>K f_2(n)$ asymptotically in $n$. For $n\in\mathbb N$ large enough and $\varepsilon'>\max\{\varepsilon,\ell_2^{-1}H(\ell_2\varepsilon),\ell_2^{-1}\mathrm{e}^2\}$, we define the events
\begin{equation}\label{eq:b1bis}
B_1^{(n)}:=\bigcup_{t=a_n}^{\lfloor n-f_1(n)\rfloor}\{S_n(t)+a_n-t\leq 0\},
\end{equation}
\begin{equation*}
B_2^{(n)}:=\bigcup_{t=\lfloor n-f_1(n)\rfloor}^{\lfloor n-\varepsilon' f_2(n)\rfloor}\{S_n(t)+a_n-t\leq 0\},
\end{equation*}
\begin{equation*}
B_3^{(n)}:=\bigcup_{t=\lfloor n-\varepsilon' f_2(n)\rfloor}^{\lfloor n-\varepsilon f_2(n)\rfloor}\{S_n(t)+a_n-t\leq 0\}
\end{equation*}
and note that
\begin{equation}\label{eq:inequalitiesbis}
P(B_3^{(n)})\leq P\left(\frac{n-A_n^*}{f_2(n)}>\varepsilon\right)\leq P(B_1^{(n)})+P(B_2^{(n)})+P(B_3^{(n)}).
\end{equation}
By Proposition \ref{ldp:supcrit} and $a_c^{(n)}/b_c^{(n)}\to +\infty$, we have
$(b_c^{(n)})^{-1}\log P(B_1^{(n)})\to-\infty$.
We shall show later on
\begin{equation}\label{eq:compasb2bis}
\limsup_{n\to\infty}\frac{1}{b_c^{(n)}}\log P(B_2^{(n)})\leq-H(\ell_2\varepsilon)
\end{equation}
and
\begin{equation}\label{eq:compasB1bis}
\lim_{n\to\infty}\frac{1}{b_c^{(n)}}\log P(B_3^{(n)})=-H(\ell_2\varepsilon).
\end{equation}
The claim then follows combining these relations with the inequality \eqref{eq:inequalitiesbis} and
the principle of the largest term (see e.g. Lemma 1.2.15 p. 7 in \cite{DZ}). We proceed by dividing the proof in three steps. Throughout
the proof we consider the quantity
\[
\kappa_n(x):=n-x f_2(n)(1+o(1)),\quad\text{$x\geq\ell_2^{-1}$, $n\in\mathbb N$.}
\]
\noindent{\it Step\,\,1:\,\,An\,\,auxiliary\,\,LDP.}\\
Let $x\geq\ell_2^{-1}$ be fixed. In this step we show that
$\{(n-a_n-S_n(\kappa_n(x)))/b_c^{(n)}\}_{n\in\mathbb N}$
obeys an LDP on $\mathbb R$ with speed $v(n):=b_c^{(n)}$ and rate function $I:=H$. Note that $n-a_n-S_n(\ell_n(x)))$,
$\ell_n(x):=n-x f_2(n)$, is distributed as $\mathrm{Bin}(n-a_n,1-\pi_n(\ell_n(x)))$. Note also that
$a_n=o(n)$, $f_2(n)\sim b_c^{(n)}=o(n)$ and by \eqref{eq:nunomenopi} (which is applicable since $b_c^{(n)}=o(p_n^{-1})$) and \eqref{eq:logblogbc1}
it follows $n(1-\pi_n(\ell_n(x)))\sim_e b_c^{(n)}$. Therefore, by Lemma \ref{ldp} we have that $\{(n-a_n-S_n(\ell_n(x)))/b_c^{(n)}\}_{n\in\mathbb N}$
obeys an LDP on $\mathbb R$ with speed $v(n):=b_c^{(n)}$ and rate function $I:=H$. Arguing as in the Step 1 of the proof of Proposition \ref{ldp:supcrit},
then the claimed LDP follows if we prove that the processes
\[
\{(n-a_n-S_n(\kappa_n(x)))/b_c^{(n)}\}_{n\in\mathbb N}\quad\text{and}\quad
\{(n-a_n-S_n(\ell_n(x)))/b_c^{(n)}\}_{n\in\mathbb N}
\]
are exponentially equivalent.
Let $\delta\in\R_+$ and $\eta\in (0,1)$ be arbitrarily fixed.
We have
\begin{align}
&\limsup_{n\to\infty}\frac{1}{b_c^{(n)}}\log P\left(|S_n(\kappa_n(x))-S_n(\ell_n(x))|>\delta b_c^{(n)}\right)\nonumber\\
&=\limsup_{n\to\infty}\frac{1}{b_c^{(n)}}\log P\left(\sum_{i\notin\mathcal{A}_n(0)}\ind\{\kappa_n(x)\wedge\ell_n(x)<Y_i^{(n)}
\leq\kappa_n(x)\vee\ell_n(x)\}>\delta b_c^{(n)}\right)\nonumber\\
&\leq\limsup_{n\to\infty}\frac{1}{b_c^{(n)}}\log P\left(\sum_{i\notin\mathcal{A}_n(0)}\ind\{(n-x f_2(n)(1+\eta)<Y_i^{(n)}
\leq n-x f_2(n)(1-\eta)\}>\delta b_c^{(n)}\right)\nonumber\\
&=\limsup_{n\to\infty}\frac{1}{b_c^{(n)}}\log P\left(\mathrm{Bin}(n-a_n,\Pi_n(x,\eta))>\delta b_c^{(n)}\right),
\label{eq:limsupbis}
\end{align}
where
\[
\Pi_n(x,\eta):=\pi_n(n-f_2(n)(1-\eta))-\pi_n(n-x f_2(n)(1+\eta)).
\]
By \eqref{eq:nunomenopi} and \eqref{eq:logblogbc1} we deduce $\frac{n\Pi_n(x,\eta)}{b_c^{(n)}}\to 0$ and so
\begin{align}
\lim_{n\to\infty}E[\mathrm{Bin}(n-a_n,\Pi_n(x,\eta))]/b_c^{(n)}=0.\label{br61}
\end{align}
Consequently, by \eqref{Penrose1.7},
for all $n$ large enough,
\begin{align}
P\left(\mathrm{Bin}(n-a_n,\Pi_n(x,\eta))
>\delta b_c^{(n)}\right)
\leq\exp\left(-\frac{\delta b_c^{(n)}}{2}\log\left(\frac{\delta b_c^{(n)}}{E[\mathrm{Bin}(n-a_n,\Pi_n(x,\eta))]}\right)\right).
\label{eq:concbis}
\end{align}
Finally, by \eqref{eq:limsupbis}, \eqref{eq:concbis} and \eqref{br61} we have
\begin{align}
\limsup_{n\to\infty}\frac{1}{b_c^{(n)}}\log P\left(|S_n(\kappa_n(x))-S_n(\ell_n(x))|>\delta b_c^{(n)}\right)=-\infty\nonumber
\end{align}
and the exponential equivalence is proved.\\
\noindent{\it Step\,\,2:\,\,Proof\,\,of\,\,\eqref{eq:compasB1bis}.}\\
We shall show later on that, for $x\geq\ell_2^{-1}$,
\begin{equation}\label{eq:exactybis}
\lim_{n\to\infty}\frac{1}{b_c^{(n)}}\log P(S_n(\kappa_n(x))\leq\lfloor\kappa_n(x)\rfloor-a_n)=-H(\ell_2 x).
\end{equation}
Since
\[
P(B_3^{(n)})\geq P(S_n(n-\varepsilon f_2(n))+a_n-\lfloor n-\varepsilon f_2(n)\rfloor\leq 0),
\]
by \eqref{eq:exactybis}
\begin{align}
\liminf_{n\to\infty}\frac{1}{b_c^{(n)}}\log P(B_3^{(n)})&\geq-H(\ell_2\varepsilon).\label{eq:lowerbis}
\end{align}
Let $t_n\in\{\lfloor n-\varepsilon' f_2(n)\rfloor,\ldots,\lfloor n-\varepsilon f_2(n)\rfloor\}$ be such that
\[
\max_{t\in\{\lfloor n-\varepsilon' f_2(n)\rfloor,\ldots,\lfloor n-\varepsilon f_2(n)\rfloor\}}P(S_n(t)\leq t-a_n)=P(S_n(t_n)\leq t_n-a_n).
\]
We have
\begin{equation}\label{eq:ineqcorrettabis}
\limsup_{n\to\infty}\frac{1}{b_c^{(n)}}\log P(S_n(t_n)\leq t_n-a_n)\leq-H(\ell_2\varepsilon).
\end{equation}
Indeed, reasoning by contradiction suppose
\[
\limsup_{n\to\infty}\frac{1}{b_c^{(n)}}\log P(S_n(t_n)\leq t_n-a_n)>-H(\ell_2\varepsilon).
\]
Letting $\{t_{n_j}\}_{j\in\mathbb N}$ denote a subsequence of $\{t_n\}_{n\in\mathbb N}$ which realizes this $\limsup$, and setting
\[
x_{n_j}:=\frac{n_j-t_{n_j}}{f_2(n_j)}
\]
we have
\[
\frac{n_j-\lfloor n_j-\varepsilon' f_2(n_j)\rfloor}{f_2(n_j)}\geq x_{n_j}\geq\frac{n_j-\lfloor n_j-\varepsilon f_2(n_j)\rfloor}{f_2(n_j)}.
\]
Therefore, by \eqref{eq:bcinfty} and \eqref{eq:intpart},
\[
\varepsilon'\geq\limsup_{j\to\infty}x_{n_j}\geq\liminf_{j\to\infty}x_{n_j}\geq\varepsilon.
\]
So, we may select a subsequence $\{x_{n_{jh}}\}_{h\in\mathbb N}\subseteq\{x_{n_j}\}_{j\in\mathbb N}$ such that $x_{n_{jh}}\to\bar x\in [\varepsilon,\varepsilon']$, as $h\to\infty$.
Consequently, for any $h\in\mathbb N$,
\[
t_{n_{jh}}=n_{jh}-\bar x f_2(n_{jh})(1+o(1))=\kappa_{n_{jh}}(\bar x).
\]
Thus by \eqref{eq:exactybis}
\[
\lim_{h\to\infty}\frac{1}{b_c^{(n_{jh})}}\log P(S_{n_{jh}}(t_{n_{jh}})\leq t_{n_{jh}}-a_{n_{jh}})=-H(\ell_2\bar x)\leq-H(\ell_2\varepsilon),
\]
where the latter inequality follows from the fact that $H$ increases on $(1,+\infty)$ and $\bar x\geq\varepsilon>\ell_2^{-1}$. This proves \eqref{eq:ineqcorrettabis}, and so
\begin{align}
&\limsup_{n\to\infty}\frac{1}{b_c^{(n)}}\log P(B_n^{(3)})\nonumber\\
&\,\,\,
\leq\limsup_{n\to\infty}\frac{1}{b_c^{(n)}}\log[(\lfloor n-\varepsilon f_2(n)\rfloor-\lfloor n-\varepsilon' f_2(n)\rfloor)P(S_n(t_n)\leq t_n-a_n)]\nonumber\\
&\,\,\,
=-H(\ell_2\varepsilon),\nonumber
\end{align}
which is the matching upper bound for \eqref{eq:lowerbis} and proves \eqref{eq:compasB1bis}. It remains to show \eqref{eq:exactybis}.
Let $\delta\in\mathbb R_+$ be arbitrarily chosen and let $n$ be so large that $o(1)<\delta$. By the definition of $\kappa_n(x)$, for all $n$ large enough,
\[
\log P(S_n(\kappa_n(x))\leq\lfloor\kappa_n(x)\rfloor-a_n)\geq\log P\left(\frac{n-a_n-S_n(\kappa_n(x))}{f_2(n)}\geq x(1+\delta)\right).
\]
It is readily checked that $\{(n-a_n-S_n(\kappa_n(x))/f_2(n)\}_{n\in\mathbb N}$
and $\{(n-a_n-S_n(\kappa_n(x)))/(\ell_2 b_c^{(n)})\}_{n\in\mathbb N}$ are exponentially equivalent.
Moreover, by the LDP of Step 1 and the Contraction Principle (see e.g. Theorem 4.2.1 p. 126 in \cite{DZ}) we have that
$\{(n-a_n-S_n(\kappa_n(x)))/(\ell_2 b_c^{(n)})\}_{n\in\mathbb N}$ obeys an LDP on $\mathbb R$ with speed $v(n):=b_c^{(n)}$
and rate function $I_2(\cdot):=H(\ell_2\cdot)$. Consequently, arguing as in the proof of Step 1 of Proposition \ref{ldp:supcrit},
$\{(n-a_n-S_n(\kappa_n(x))/f_2(n)\}_{n\in\mathbb N}$ obeys an LDP on $\mathbb R$ with speed $v(n):=b_c^{(n)}$
and rate function $I_2$. Therefore
\begin{align}
&\liminf_{n\to\infty}\frac{1}{b_c^{(n)}}\log P(S_n(\kappa_n(x))\leq\lfloor\kappa_n(x)\rfloor-a_n)\nonumber\\
&\,\,\,\,\,\,
\geq\liminf_{n\to\infty}\frac{1}{b_c^{(n)}}\log P\left(\frac{n-a_n-S_n(\kappa_n(x))}{f_2(n)}>x(1+\delta)\right)\nonumber\\
&\,\,\,\,\,\,
\geq-\inf_{y\in (x(1+\delta),+\infty)}H(\ell_2 y)
=-H(\ell_2 x(1+\delta))\label{eq:20bisbis},
\end{align}
where for the equality in \eqref{eq:20bisbis} we used that $H$ is continuously increasing on $(1,+\infty)$. Taking the supremum over $\delta>0$,
we deduce the lower bound
\begin{align}
\liminf_{n\to\infty}\frac{1}{b_c^{(n)}}\log P(S_n(\kappa_n(x))\leq\lfloor\kappa_n(x)\rfloor-a_n)&\geq-H(\ell_2 x).\nonumber
\end{align}
Since $H(1)=0$ the matching upper bound with $x=\ell_2^{-1}$ is trivially true. It remains to prove the matching upper bound for $x>\ell_2^{-1}$.
Take $\delta\in (0,1-(\ell_2 x)^{-1})$ and let $n$ be so large that $-\varepsilon<o(1)$. Using again the LDP for $\{(n-a_n-S_n(\kappa_n(x))/f_2(n)\}_{n\in\mathbb N}$, we have
\begin{align}
&\limsup_{n\to\infty}\frac{1}{b_c^{(n)}}\log P(S_n(\kappa_n(x))\leq\lfloor\kappa_n(x)\rfloor-a_n)\nonumber\\
&\,\,\,\,\,\,
\leq\limsup_{n\to\infty}\frac{1}{b_c^{(n)}}\log P\left(\frac{n-a_n-S_n(\kappa_n(x))}{f_2(n)}\geq x(1-\delta)\right)\nonumber\\
&\,\,\,\,\,\,
\leq-\inf_{[x(1-\delta),+\infty)}H(\ell_2 y)=-H(\ell_2 x(1-\delta))\nonumber
\end{align}
and the matching upper bound follows by letting $\delta$ tend to zero.\\
\noindent{\it Step\,\,3:\,\,Proof\,\,of\,\,\eqref{eq:compasb2bis}.}\\
Letting $S_n'$ denote the process defined by \eqref{eq:snprime},
we have
\begin{align}
P(B_2^{(n)})&\leq P(S_n(n-f_1(n))\leq\lfloor n-\varepsilon' f_2(n)\rfloor-a_n)\nonumber\\
&\leq P(S_n'(n-f_1(n))\leq\lfloor n-\varepsilon' f_2(n)\rfloor)\nonumber\\
&=P(\mathrm{Bin}(n,1-\pi_n(n-f_1(n)))\geq n-\lfloor n-\varepsilon' f_2(n)\rfloor)\nonumber\\
&=P(\mathrm{Bin}(n,1-\pi_n(n-f_1(n)))\geq\lceil\varepsilon' f_2(n)\rceil),\nonumber
\end{align}
where the latter relation follows by \eqref{eq:intpart}.
By \eqref{eq:nunomenopi}, \eqref{eq:logblogbc1} and the definition of $\varepsilon'$ one has
\[
\lim_{n\to\infty}\frac{\lceil\varepsilon' f_2(n)\rceil}{n(1-\pi_n(n-f_1(n)))}
=\ell_2\varepsilon'>\mathrm{e}^2.
\]
Therefore, by \eqref{Penrose1.7},
for all $n$ large enough,
\begin{align}
P(B_2^{(n)})&\leq\exp\left(-\frac{\lceil\varepsilon' f_2(n)\rceil}{2}\log\left(\frac{\lceil\varepsilon' f_2(n)\rceil}{n(1-\pi_n(n-f_1(n)))}\right)\right)\nonumber\\
&\leq\exp\left(-\lceil\varepsilon' f_2(n)\rceil\right)\leq\exp(-\varepsilon' f_2(n)).\nonumber
\end{align}
By this inequality we easily have \eqref{eq:compasb2bis} since $\varepsilon'>\ell_2^{-1}H(\ell_2\varepsilon)$.\\
\noindent{\it Proof\,\,of\,\,\eqref{eq:lim2}.}\\
Let $\varepsilon\in (0,\ell_2^{-1})$ and $\varepsilon'>\ell_2^{-1}$ be such that $H(\ell_2\varepsilon')>H(\ell_2^{-1}\varepsilon)$. We have
\begin{align}
P\left(\frac{n-A_n^*}{f_2(n)}\leq\varepsilon\right)
&=P(A_n(t)>t,\quad\forall t=a_n,\ldots,\lfloor n-\varepsilon f_2(n)\rfloor)\nonumber\\
&=P(B_4^{(n)}\cap B_5^{(n)}),\nonumber
\end{align}
where
\[
B_4^{(n)}:=\bigcap_{t=a_n}^{\lfloor n-\varepsilon' f_2(n)\rfloor}\{A_n(t)>t\}\quad\text{and}\quad
B_5^{(n)}:=\bigcap_{t=\lfloor n-\varepsilon' f_2(n)\rfloor}^{\lfloor n-\varepsilon f_2(n)\rfloor}\{A_n(t)>t\}.
\]
We proceed by diving the proof in two steps. In the first step we prove
\begin{equation}\label{eq:pdn}
\lim_{n\to\infty}\frac{1}{b_c^{(n)}}\log P(B_5^{(n)})=-H(\ell_2\varepsilon)
\end{equation}
and in the second step we conclude the proof of \eqref{eq:lim2}.\\
\noindent{\it Step\,\,1:\,\,Proof\,\,of\,\,\eqref{eq:pdn}.}\\
Since
\[
P(S_n(n-\varepsilon' f_2(n))+a_n>\lfloor n-\varepsilon f_2(n)\rfloor)
\leq P(B_5^{(n)})\leq P(S_n(n-\varepsilon f_2(n))+a_n>\lfloor n-\varepsilon f_2(n)\rfloor),
\]
the claim follows if we prove that, for any $x\in\mathbb R_+$,
\begin{equation}\label{eq:pdn1}
\lim_{n\to\infty}\frac{1}{b_c^{(n)}}\log P(S_n(n-x f_2(n))+a_n>\lfloor n-\varepsilon f_2(n)\rfloor)=-H(\ell_2\varepsilon).
\end{equation}
By \eqref{eq:intpart} we have
\begin{align}
P(S_n(n-x f_2(n))+a_n>\lfloor n-\varepsilon f_2(n)\rfloor)&=P(n-a_n-S_n(n-x f_2(n))<\lceil\varepsilon f_2(n)\rceil)\nonumber\\
&=P\left(\frac{n-a_n-S_n(n-x f_2(n))}{f_2(n)}\leq\varepsilon\right).\nonumber
\end{align}
The claim \eqref{eq:pdn1} follows by this relation
noticing that arguing as in the Step 2 of the proof of \eqref{eq:lim1} one has that
the process $\Biggl\{\frac{n-a_n-S_n(n-x f_2(n))}{f_2(n)}\Biggr\}_{n\in\mathbb N}$ obeys an LDP on $\mathbb R$ with speed
$v(n):=b_c^{(n)}$ and rate function $I_2$.\\
\noindent{\it Step\,\,2:\,\,Conclusion\,\,of\,\,the\,\,proof\,\,of\,\,\eqref{eq:lim2}.}\\
By \eqref{eq:lim1}
\[
\lim_{n\to\infty}\frac{1}{b_c^{(n)}}\log(1-P(B_4^{(n)}))=-H(\ell_2\varepsilon').
\]
Combining this with \eqref{eq:pdn}, since $\varepsilon'>\ell_2^{-1}$ is chosen in such a way that $H(\ell_2\varepsilon')>H(\ell_2\varepsilon)$, we have $(1-P(B_4^{(n)}))/P(B_5^{(n)})\to 0$. Since
\[
\frac{P(B_4^{(n)})+P(B_5^{(n)})-1}{P(B_5^{(n)})}=1-\frac{1-P(B_4^{(n)})}{P(B_5^{(n)})}\leq\frac{P(B_4^{(n)}\cap B_5^{(n)})}{P(B_5^{(n)})}\leq 1
\]
we deduce $P(B_4^{(n)}\cap B_5^{(n)})\sim_e P(B_5^{(n)})$ and so
\begin{align}
\lim_{n\to\infty}\frac{1}{b_c^{(n)}}\log P\left(\frac{n-A_n^*}{f_2(n)}\leq\varepsilon\right)&=\lim_{n\to\infty}\frac{1}{b_c^{(n)}}\log P(B_4^{(n)}\cap B_5^{(n)})\nonumber\\
&=\lim_{n\to\infty}\frac{1}{b_c^{(n)}}\log P(B_5^{(n)})=-H(\ell_2\varepsilon),\nonumber
\end{align}
which concludes the proof.

\subsubsection{Proof of Theorem \ref{ldp2}}

We divide the proof in two steps. In the first step we prove the large deviation lower bound
and in the second step we prove the large deviation upper bound.\\
\noindent{\it Step\,\,1:\,\,large\,\,deviation\,\,lower\,\,bound.}\\
Let $O \subseteq\R$ be an open set. If $O\subseteq\mathbb R_{-}$, then the claim is obvious
and so we assume $O\cap [0,\infty)\neq\emptyset$. By \eqref{eq:ratiopos}
we have
\[
P\left(\frac{n-A_n^*}{f_2(n)}\in O\right)=P\left(\frac{n-A_n^*}{f_2(n)}\in O\cap [0,\infty)\right).
\]
If $\ell_2^{-1}\in O$, then $\ell_2^{-1}\in O\cap\R_+$ and since $O\cap\R_+$ is open there exist $\varepsilon_1<\ell_2^{-1}<\varepsilon_2$ such that $(\varepsilon_1,\varepsilon_2]\subset O\cap\R_+$.
So
\begin{align}
P\left(\frac{n-A_n^*}{f_2(n)}\in O\right)&\geq P\left(\frac{n-A_n^*}{f_2(n)}\in (\varepsilon_1,\varepsilon_2]\right)\nonumber\\
&=P\left(\frac{n-A_n^*}{f_2(n)}\leq\varepsilon_2\right)-P\left(\frac{n-A_n^*}{f_2(n)}\leq\varepsilon_1\right).\label{eq:disinbis}
\end{align}
By Proposition \ref{ldp:supcritbis} we have that, for $\eta\in (0,H(\ell_2\varepsilon_1)\wedge H(\ell_2\varepsilon_2))$ arbitrarily fixed, there exists $n_\eta$ such that for any $n>n_\eta$
\begin{equation}\label{eq:exp1terbis}
1-\mathrm{e}^{-(H(\ell_2\varepsilon_2)-\eta)b_c^{(n)}}<P\left(\frac{n-A_n^*}{f_2(n)}\leq\varepsilon_2\right)<1-\mathrm{e}^{-(H(\ell_2\varepsilon_2)+\eta)b_c^{(n)}}
\end{equation}
and
\begin{equation}\label{eq:exp1bis}
\mathrm{e}^{-(H(\ell_2\varepsilon_1)+\eta)b_c^{(n)}}<P\left(\frac{n-A_n^*}{f_2(n)}\leq\varepsilon_1\right)<\mathrm{e}^{-(H(\ell_2\varepsilon_1)-\eta)b_c^{(n)}}.
\end{equation}
By \eqref{eq:disinbis}, \eqref{eq:exp1terbis} and \eqref{eq:exp1bis}
we easily have
\[
\liminf_{n\to\infty}\frac{1}{b_c^{(n)}}\log P\left(\frac{n-A_n^*}{f_2(n)}\in O\right)\geq 0=-\inf_{y\in O} I_2(y),
\]
and the large deviation lower bound for the case $\ell_2^{-1}\in O$ is proved. If $\ell_2^{-1}\notin O$, then $O\cap [0,+\infty)=(O\cap [0,\ell_2^{-1}))\cup(O\cap(\ell_2^{-1},+\infty))$.
Since $O\cap [0,+\infty)\neq\emptyset$ then $O\cap [0,\ell_2^{-1})\neq\emptyset$
and/or $O\cap(\ell_2^{-1},+\infty)\neq\emptyset$. In particular note that if $O\cap [0,\ell_2^{-1})\neq\emptyset$ then $O\cap (0,\ell_2^{-1})\neq\emptyset$, since $O$
is open. Thus, $O\cap (0,\ell_2^{-1})\neq\emptyset$ and/or $O\cap(\ell_2^{-1},+\infty)\neq\emptyset$. In the following, we suppose $O\cap (0,\ell_2^{-1})\neq\emptyset$ and
$O\cap(\ell_2^{-1},+\infty)\neq\emptyset$, however the proof can be easily adapted to the case when one of these two intersections is empty.
Let $x\in O\cap (0,\ell_2^{-1})$ and $z\in O\cap(\ell_2^{-1},+\infty)$. For all $\delta\in (0,x\wedge(\ell_2^{-1}-x)\wedge(z-\ell_2^{-1}))$ small enough,
we have $O\cap (0,\ell_2^{-1})\supset (x-\delta,x+\delta]$ and $O\cap(\ell_2^{-1},+\infty)\supset (z-\delta,z+\delta]$.
Therefore,
\begin{align}
P\left(\frac{n-A_n^*}{f_2(n)}\in O\right)&\geq P\left(\frac{n-A_n^*}{f_2(n)}\in O\cap (0,\ell_2^{-1})\right)+P\left(\frac{n-A_n^*}{f_2(n)}\in O\cap (\ell_2^{-1},+\infty)\right)\nonumber\\
&\geq P\left(\frac{n-A_n^*}{f_2(n)}\leq x+\delta\right)-P\left(\frac{n-A_n^*}{f_2(n)}\leq x-\delta\right)\nonumber\\
&\,\,\,\,\,\,
\,\,\,\,\,\,
+P\left(\frac{n-A_n^*}{f_2(n)}\leq z+\delta\right)-P\left(\frac{n-A_n^*}{f_2(n)}\leq z-\delta\right).\label{eq:lb2}
\end{align}
By Proposition \ref{ldp:supcritbis}, we have that, for any $\eta\in (0,2^{-1}[H(\ell_2(x-\delta))-H(\ell_2(x+\delta))]\wedge
2^{-1}[H(\ell_2(z+\delta))-H(\ell_2(z-\delta))])$, there exists $n_\eta$ such that for all $n>n_\eta$
the inequality \eqref{eq:exp1terbis} holds with $\varepsilon_2=z-\delta$ and $\varepsilon_2=z+\delta$
and the inequality \eqref{eq:exp1bis} holds with $\varepsilon_1=x-\delta$ and $\varepsilon_1=x+\delta$.
So, for all $n$ large enough,
\begin{align}
P\left(\frac{n-A_n^*}{f_2(n)}\leq z+\delta\right)&-P\left(\frac{n-A_n^*}{f_2(n)}\leq z-\delta\right)
\geq\mathrm{e}^{-(H(\ell_2(z-\delta))+\eta)b_c^{(n)}}-\mathrm{e}^{-(H(\ell_2(z+\delta))-\eta)b_c^{(n)}}\nonumber\\
&=\mathrm{e}^{-(H(\ell_2(z-\delta))+\eta)b_c^{(n)}}\left(1-\mathrm{e}^{-(H(\ell_2(z+\delta))-H(\ell_2(z-\delta))-2\eta)b_c^{(n)}}\right)\nonumber
\end{align}
and
\begin{align}
P\left(\frac{n-A_n^*}{f_2(n)}\leq x+\delta\right)&-P\left(\frac{n-A_n^*}{f_2(n)}\leq x-\delta\right)
\geq\mathrm{e}^{-(H(\ell_2(x+\delta))+\eta)b_c^{(n)}}-\mathrm{e}^{-(H(\ell_2(x-\delta))-\eta)b_c^{(n)}}\nonumber\\
&=\mathrm{e}^{-(H(\ell_2(x+\delta))+\eta)b_c^{(n)}}\left(1-\mathrm{e}^{-(H(\ell_2(x-\delta))-H(\ell_2(x+\delta))-2\eta)b_c^{(n)}}\right).\nonumber
\end{align}
By these inequalities and \eqref{eq:lb2} we have
\begin{align}
&\liminf_{n\to\infty}\frac{1}{b_c^{(n)}}\log P\left(\frac{n-A_n^*}{f_2(n)}\in O\right)\nonumber\\
&\,\,\,
\geq\max\Biggl\{\liminf_{n\to\infty}\frac{1}{b_c^{(n)}}\log\left(P\left(\frac{n-A_n^*}{f_2(n)}\leq x+\delta\right)-P\left(\frac{n-A_n^*}{f_2(n)}\leq x-\delta\right)\right),\nonumber\\
&\,\,\,\,\,\,
\,\,\,\,\,\,
\,\,\,\,\,\,
\,\,\,\,\,\,
\,\,\,\,\,\,
\liminf_{n\to\infty}\frac{1}{b_c^{(n)}}\log \left(P\left(\frac{n-A_n^*}{f_2(n)}\leq z+\delta\right)-P\left(\frac{n-A_n^*}{f_2(n)}\leq z-\delta\right)\right)\Biggr\}\nonumber\\
&\geq(-H(\ell_2(x+\delta))-\eta)\vee(-H(\ell_2(z-\delta))-\eta).\nonumber
\end{align}
So, letting first $\eta$ and $\delta$ tend to zero and taking then the supremum over $\{x\in O\cap (0,\ell_2^{-1}), z\in O\cap (\ell_2^{-1},+\infty)\}$, we have
\begin{align}
\liminf_{n\to\infty}\frac{1}{b_c^{(n)}}\log P\left(\frac{n-A_n^*}{f_2(n)}\in O\right)&\geq-\inf_{x\in O\cap (0,\ell_2^{-1})}H(\ell_2 x)\wedge\inf_{z\in O\cap (\ell_2^{-1},+\infty)}H(\ell_2z)\nonumber\\
&=-\inf_{y\in O\cap [0,\infty)}H(\ell_2 y)=-\inf_{y\in O}I_2(y),\nonumber
\end{align}
and the large deviation lower bound is proved.\\
\noindent{\it Step\,\,2:\,\,large\,\,deviation\,\,upper\,\,bound.}\\
Let $C\subseteq\mathbb R$ be a closed set. If $C\subset\mathbb R_{-}$, then by \eqref{eq:ratiopos} the large deviation upper bound is trivial. Therefore, we assume $C\cap [0,\infty)\neq\emptyset$.
If $\ell_2^{-1}\in C$, then the large deviation upper bound is again trivial. If $\ell_2^{-1}\notin C$, then $C\cap [0,+\infty)=(C\cap [0,\ell_2^{-1}))\cup (C\cap(\ell_2^{-1},+\infty))$.
In the following, we suppose $C\cap [0,\ell_2^{-1})\neq\emptyset$ and $C\cap(\ell_2^{-1},+\infty)\neq\emptyset$, however the proof can be easily adapted to the case
when one of these two intersections is empty. Set $M:=\max\{x:\,\,x\in C\cap [0,\ell_2^{-1})\}$ and $m:=\min\{x:\,\,x\in C\cap (\ell_2^{-1},+\infty)\}$. Then $C\cap [0,\ell_2^{-1})\subset [0,M]$ and
$C\cap(\ell_2^{-1},+\infty)\subset (m-\delta,+\infty)$, for all $\delta\in\mathbb R_+$ small enough. Consequently,
\begin{align}
P\left(\frac{n-A_n^*}{f_2(n)}\in C\right)&\leq P\left(\frac{n-A_n^*}{f_2(n)}\leq M\right)+P\left(\frac{n-A_n^*}{f_2(n)}>m-\delta\right)\nonumber
\end{align}
and so by Proposition \ref{ldp:supcritbis} and the principle of the largest term (see e.g. Lemma 1.2.15 p. 7 in \cite{DZ})
\begin{align}
\limsup_{n\to\infty}\frac{1}{b_c^{(n)}}\log P\left(\frac{n-A_n^*}{f_2(n)}\in C\right)&\leq(-H(\ell_2 M))\vee (-H(\ell_2(m-\delta))).\nonumber
\end{align}
Therefore, letting $\delta$ tend to zero,
\begin{align}
\limsup_{n\to\infty}\frac{1}{b_c^{(n)}}\log P\left(\frac{n-A_n^*}{f_2(n)}\in C\right)&\leq (-H(\ell_2 M))\vee(-H(\ell_2 m))\nonumber\\
&=-\inf_{y\in C\cap [0,\ell_2^{-1})}I_2(y)\vee -\inf_{y\in C\cap (\ell_2^{-1},\infty)}I_2(y)\nonumber\\
&=-\inf_{y\in C}I_2(y),\nonumber
\end{align}
and the large deviation upper bound is proved.

\subsection{Proofs of Proposition \ref{ldp:supcritpenta} and Theorem \ref{ldp5}}

\subsubsection{Proof of Proposition \ref{ldp:supcritpenta}}

Let $f_1(n)$ be the function defined in the statement of Theorem \ref{ldp1}. Since for $n$ large enough
$f_1(n)=(g(n)a_c^{(n)})/(n p_n)$, by the second relation in \eqref{eq:f3property} $f_1(n)>f_3(n)$ for all $n$ large enough. Hereafter, we chose $g$ in such a way that
$f_1(n)=o(p_n^{-1})$. For $n\in\mathbb N$ large enough, let $B_i^{(n)}$, $i=1,2,3$, be defined by
\eqref{eq:nuovoB1}, \eqref{eq:nuovoB2}, \eqref{eq:b3} with $\varepsilon=1$, and set
\[
B_4^{(n)}:=\bigcup_{t=\lfloor n-f_1(n)\rfloor}^{\lfloor n-\varepsilon f_3(n)\rfloor}\{S_n(t)+a_n-t\leq 0\}.
\]
We have
\begin{equation}\label{eq:uplowpanomega}
P(B_4^{(n)})\leq P\left(\frac{n-A_n^*}{f_3(n)}>\varepsilon\right)\leq P(B_1^{(n)})+P(B_2^{(n)})+P(B_3^{(n)})+P(B_4^{(n)}).
\end{equation}
Since
\[
P(S_n(n-\varepsilon f_3(n))+a_n\leq \lfloor n-\varepsilon f_3(n)\rfloor)\leq P(B_4^{(n)})\leq P(S_n(n-f_1(n))+a_n\leq \lfloor n-\varepsilon f_3(n)\rfloor),
\]
by \eqref{eq:intpart} we have
\begin{equation}\label{eq:uplowpb4}
P(n-a_n-S_n(n-\varepsilon f_3(n))\geq\lceil\varepsilon f_3(n)\rceil)\leq P(B_4^{(n)})\leq P(n-a_n-S_n(n-f_1(n))\geq\lceil\varepsilon f_3(n)\rceil).
\end{equation}
Let $h$ be such that $h(n)=o(p_n^{-1})$. By \eqref{eq:nunomenopi}
\begin{equation}\label{eq:rate1menopi}
n(1-\pi_n(n-h(n)))\sim_e b_c^{(n)'}.
\end{equation}
By \eqref{eq:logblogbc1} we have $b_c^{(n)}\sim_e b_c^{(n)'}$, and so
$b_c^{(n)'}\to b$.
Note that $f_3(n)/b_c^{(n)'}\to\infty$ by the first relation in \eqref{eq:f3property},
and $f_3(n)=o(n)$ since $a_c^{(n)}=o(n)$ and $f_3(n)=o(a_c^{(n)})$ by the second relation in \eqref{eq:f3property}. So by \eqref{eq:nuovorate} and \eqref{eq:rate1menopi}
\begin{align}
&\log P(n-a_n-S_n(n-h(n))\geq\lceil\varepsilon f_3(n)\rceil)=\log P(\mathrm{Bin}(n-a_n,1-\pi_n(n-h(n)))\geq\lceil\varepsilon f_3(n)\rceil)\nonumber\\
&\,\,\,\,\,\,
\,\,\,\,\,\,
\,\,\,\,\,\,
\sim_e\lceil\varepsilon f_3(n)\rceil\log\left(\frac{(n-a_n)(1-\pi_n(n-h(n)))}{f_3(n)}\right)\nonumber\\
&\,\,\,\,\,\,
\,\,\,\,\,\,
\,\,\,\,\,\,
\sim_e\varepsilon f_3(n)\log\left(\frac{b_c^{(n)'}}{f_3(n)}\right).\label{eq:Stir1}
\end{align}
Since $f_1(n)$ and $f_3(n)$ are both $o(p_n^{-1})$,
by \eqref{eq:uplowpb4}
\begin{equation}\label{eq:PB4}
\lim_{n\to\infty}\frac{1}{-f_3(n)\log(b_c^{(n)}/f_3(n))}\log P(B_n^{(4)})=-\varepsilon.
\end{equation}
The claim follows by the inequality \eqref{eq:uplowpanomega}, the principle of the largest term
(see e.g. Lemma 1.2.15 p. 7 in \cite{DZ}), relations
\eqref{eq:compasB1}, \eqref{eq:compasb2}, \eqref{eq:compasb3}, \eqref{eq:PB4} and the fact that
$-f_3(n)\log(b_c^{(n)}/f_3(n))=o(a_c^{(n)})$ (this latter relation easily follows by the first relation in \eqref{eq:hyp}, \eqref{eq:tc},
the definition of $b_c^{(n)}$ and the second relation in \eqref{eq:f3property}).

\subsubsection{Proof of Theorem \ref{ldp5}}

We divide the proof in two steps. In the first step we prove the large deviation lower bound
and in the second step we prove the large deviation upper bound.\\
\noindent{\it Step\,\,1:\,\,large\,\,deviation\,\,lower\,\,bound.}\\
Let $O\subseteq\R$ be an open set. If $O\subseteq\mathbb R_{-}$, then the claim is obvious
and so we assume $O\cap [0,\infty)\neq\emptyset$. By \eqref{eq:ratiopos}
we have
\[
P\left(\frac{n-A_n^*}{f_3(n)}\in O\right)=P\left(\frac{n-A_n^*}{f_3(n)}\in O\cap [0,\infty)\right).
\]
If $0\in O$, then there exists $\delta\in\mathbb R_+$ such that $(-\delta,\delta)\subset O$ and for a fixed $0<\varepsilon<\delta$,
\begin{align*}
P\left(\frac{n-A_n^*}{f_3(n)}\in O\right)\geq P\left(\frac{n-A_n^*}{f_3(n)}\in [0,\varepsilon]\right).
\end{align*}
By Proposition \ref{ldp:supcritpenta}, for any $\eta\in (0,\varepsilon)$, there exists $n_\eta$ such that for any $n>n_\eta$
\begin{equation}
\label{eq:exp1ter}
\mathrm{e}^{(\varepsilon+\eta)f_3(n)\log(b_c^{(n)}/f_3(n))}<P\left(\frac{n-A_n^*}{f_3(n)}>\varepsilon\right)<\mathrm{e}^{(\varepsilon-\eta)f_3(n)\log(b_c^{(n)}/f_3(n))}.
\end{equation}
Therefore
\begin{align}
&\liminf_{n\to\infty}\frac{1}{-f_3(n)\log(b_c^{(n)}/f_3(n))}\log P\left(\frac{n-A_n^*}{f_3(n)}\in O\right)\nonumber\\
&\,\,\,\,\,\,
\geq\lim_{n\to\infty}\frac{1}{-f_3(n)\log(b_c^{(n)}/f_3(n))}\log\left(1-P\left(\frac{n-A_n^*}{f_3(n)}>\varepsilon\right)\right)
=0=-\inf_{x\in O}I_3(x),\nonumber
\end{align}
and the large deviation lower bound for the case $0\in O$ is proved. If $0\notin O$, then $O\cap [0,+\infty)=O\cap\R_+$. Let $x\in O\cap\R_+$
be arbitrarily fixed. Since $O\cap\mathbb R_+$ is open, there exists
$\delta\in (0,x)$
such that $(x-\delta,x+\delta]\subset O\cap\R_+$. Therefore,
\begin{align}
P\left(\frac{n-A_n^*}{f_3(n)}\in O\right)&\geq P\left(\frac{n-A_n^*}{f_3(n)}>x-\delta\right)-P\left(\frac{n-A_n^*}{f_3(n)}>x+\delta\right).\nonumber
\end{align}
Note that for any $\eta\in (0,\delta)$, there exists $n_\eta$ such that for any $n>n_\eta$ \eqref{eq:exp1ter} holds with $x\pm\delta$ in place of $\varepsilon$, and so
\begin{align}
P\left(\frac{n-A_n^*}{f_3(n)}\in O\right)&\geq\mathrm{e}^{(\eta+x-\delta)f_3(n)\log(b_c^{(n)}/f_3(n))}-\mathrm{e}^{(x+\delta-\eta)f_3(n)\log(b_c^{(n)}/f_3(n))}\nonumber\\
&=\mathrm{e}^{(x-\delta+\eta)f_3(n)\log(b_c^{(n)}/f_3(n))}\left(1-\mathrm{e}^{2(\delta-\eta)f_3(n)\log(b_c^{(n)}/f_3(n))}\right)\nonumber
,\quad\text{for any $n>n_\eta$.}\nonumber
\end{align}
Therefore, taking the logarithm on this inequality, dividing then by $-f_3(n)\log(b_c^{(n)}/f_3(n))$, and letting first $n$ tend to $\infty$ and second $\eta$ tend to zero,
\[
\liminf_{n\to\infty}\frac{1}{-f_3(n)\log(b_c^{(n)}/f_3(n))}\log P\left(\frac{n-A_n^*}{f_3(n)}\in O\right)\geq -x+\delta\geq-x.
\]
The large deviation lower bound follows taking the supremum over all $x\in O\cap\mathbb R_+$ on this relation.\\
\noindent{\it Step\,\,2:\,\,large\,\,deviation\,\,upper\,\,bound.}\\
Let $C\subseteq\mathbb R$ be a closed set. If $C\subset\mathbb R_{-}$,
then by \eqref{eq:ratiopos} the large deviation upper bound is trivial. Therefore, we assume $C\cap [0,\infty)\neq\emptyset$.
If $0\in C$, then the large deviation upper bound is again trivial. If $0\notin C$, then $C\cap [0,+\infty)=C\cap\R_+$.
Let $m:=\min\{x:\,\,x\in C\cap\R_+\}$. For all $\delta\in (0,m)$, we have $C\cap [0,+\infty)\subset (m-\delta,+\infty)$ and so
by Proposition \ref{ldp:supcritpenta} we deduce
\begin{align}
\limsup_{n\to\infty}\frac{1}{-f_3(n)\log(b_c^{(n)}/f_3(n))}\log P\left(\frac{n-A_n^*}{f_3(n)}\in C\right)\leq-m+\delta.\nonumber
\end{align}
Letting $\delta$ tend to zero, we have
\begin{align}
\limsup_{n\to\infty}\frac{1}{-f_3(n)\log(b_c^{(n)}/f_3(n))}\log P\left(\frac{n-A_n^*}{f_3(n)}\in C\right)\leq-m=-\inf_{x\in C}I_3(x),\nonumber
\end{align}
and the large deviation upper bound is proved.

\subsection{Proofs of Proposition \ref{ldp:supcritter} and Theorem \ref{ldp3}}

\subsubsection{Proof of Proposition \ref{ldp:supcritter}}

\noindent{\it Proof\,\,of\,\,$(i)$.}\\
Let $f_1$ be the function defined in the statement of Theorem \ref{ldp1} and let $\varepsilon\in\R_+$ be
arbitrarily fixed. By \eqref{acpninfty}
we have $f_1(n)=\frac{g(n)a_c^{(n)}}{n p_n}$ for $n$ large enough, and so $f_4(n)\ll f_1(n)$.
We choose $g$ in such a way that $f_1(n)=o(p_n^{-1})$ (and so, in particular, $f_1(n)/n\to 0$).
For $n\in\mathbb N$ sufficiently large, let $B_1^{(n)}$ be the event defined by \eqref{eq:b1bis} and set
\begin{equation}\label{eq:b2ter}
B_2^{(n)}:=\bigcup_{t=\lfloor n-f_1(n)\rfloor}^{\lfloor n-\varepsilon f_4(n)\rfloor}\{S_n(t)+a_n-t\leq 0\}.
\end{equation}
Clearly,
\begin{equation}\label{eq:ineqter}
P(B_2^{(n)})\leq P\left(\frac{n-A_n^*}{f_4(n)}>\varepsilon\right)\leq P(B_1^{(n)})+P(B_2^{(n)}).
\end{equation}
In the next steps, we shall show
\begin{equation}\label{eq:compasB1ter}
\lim_{n\to\infty}\frac{1}{-\log b_c^{(n)}}\log P(B_1^{(n)})=-\infty
\end{equation}
and
\begin{equation}\label{eq:compasb2ter}
\lim_{n\to\infty}\frac{1}{-\log b_c^{(n)}}\log P(B_2^{(n)})=-\lceil\ell_4\varepsilon\rceil.
\end{equation}
The claim then follows combining these relations with \eqref{eq:ineqter} and the principle of the largest term (see e.g. Lemma 1.2.15 p. 7 in \cite{DZ}).\\
\noindent{\it Step\,\,1:\,\,Proof\,\,of\,\,\eqref{eq:compasB1ter}.}\\
By the second relation in \eqref{eq:logblogbc}, $-\log b_c^{(n)}\sim n p_n$ and so by \eqref{acpninfty} we deduce
$a_c^{(n)}/-\log b_c^{(n)}\to+\infty$. Combining this with Proposition \ref{ldp:supcrit} and \eqref{eq:pdopo}
\begin{align}
\lim_{n\to\infty}\frac{1}{-\log b_c^{(n)}}\log P(B_1^{(n)})&=
\lim_{n\to\infty}\frac{a_c^{(n)}}{-\log b_c^{(n)}}\lim_{n\to\infty}\frac{1}{a_c^{(n)}}\log P(B_1^{(n)})\nonumber\\
&=-J(x_0)\lim_{n\to\infty}\frac{a_c^{(n)}}{-\log b_c^{(n)}}=-\infty.\nonumber
\end{align}
\noindent{\it Step\,\,2:\,\,Proof\,\,of\,\,\eqref{eq:compasb2ter}.}\\
Note that
\[
P(S_n(n-\varepsilon f_4(n))+a_n\leq\lfloor n-\varepsilon f_4(n)\rfloor)
\leq P(B_2^{(n)})\leq P(S_n(n-f_1(n))+a_n\leq\lfloor n-\varepsilon f_4(n)\rfloor)
\]
and so by \eqref{eq:intpart}
\begin{equation}
\label{eq:292}
P(n-a_n-S_n(n-\varepsilon f_4(n))\geq\lceil\varepsilon f_4(n)\rceil)
\leq P(B_2^{(n)})\leq P(n-a_n-S_n(n-f_1(n))\geq\lceil\varepsilon f_4(n)\rceil).
\end{equation}
The claim follows if we check that, for $h(n)=o(p_n^{-1})$,
\begin{equation}\label{eq:limfinnuovo}
\lim_{n\to\infty}\frac{1}{-\log b_c^{(n)}}\log P(n-a_n-S_n(n-h(n))\geq\lceil\varepsilon f_4(n)\rceil)=-\lceil\ell_4\varepsilon\rceil.
\end{equation}
Since $b_c^{(n)}\to 0$, by \eqref{eq:logblogbc1} we have $b_c^{(n)'}\to 0$ and so by
\eqref{eq:nunomenopi} we deduce $(n-a_n)(1-\pi_n(n-h(n)))\to 0$. Therefore by \eqref{eq:nuovorate}
\begin{align}
&\log P(n-a_n-S_n(n-h(n))\geq\lceil\varepsilon f_4(n)\rceil)=\log P(\mathrm{Bin}(n-a_n,1-\pi_n(n-h(n)))\geq\lceil\varepsilon f_4(n)\rceil)\nonumber\\
&\,\,\,\,\,\,
\sim_e\lceil\varepsilon f_4(n)\rceil\log b_c^{(n)'}.\label{eq:rellog}
\end{align}
Relation \eqref{eq:limfinnuovo} follows by \eqref{eq:rellog} and the first relation in \eqref{eq:logblogbc}.\\
\noindent{\it Proof\,\,of\,\,$(ii)$.}\\
Let $f_1$ be the function defined in the statement of Theorem \ref{ldp1} and let $\varepsilon\in\mathbb R_+$ be arbitrarily fixed.
Here again,
$f_1(n)=\frac{g(n)a_c^{(n)}}{n p_n}$ for $n$ large enough and we chose $g$ in such a way that $f_1(n)=o(p_n^{-1})$. Note that by \eqref{eq:f4o} one has
$f_4(n)\ll f_1(n)$. For $n\in\mathbb N$ sufficiently large, let $B_1^{(n)}$ be the event defined by \eqref{eq:b1bis} and let $B_2^{(n)}$ be the event
defined by \eqref{eq:b2ter}. Clearly, we still have \eqref{eq:ineqter}.
In the next steps, we shall show
\begin{equation}\label{eq:compasB1terBIS}
\lim_{n\to\infty}\frac{1}{-f_4(n)\log(b_c^{(n)}/f_4(n))}\log P(B_1^{(n)})=-\infty
\end{equation}
and
\begin{equation}\label{eq:compasb2terbis}
\lim_{n\to\infty}\frac{1}{-f_4(n)\log(b_c^{(n)}/f_4(n))}\log P(B_2^{(n)})=-\varepsilon.
\end{equation}
The claim then follows combining these relations with \eqref{eq:ineqter} and the principle of the largest term.
\noindent{\it Step\,\,1:\,\,Proof\,\,of\,\,\eqref{eq:compasB1terBIS}.}\\
Note that
$-f_4(n)\log(b_c^{(n)}/f_4(n))=o(a_c^{(n)})$ (this easily follows by the first relation in \eqref{eq:hyp},
\eqref{eq:tc},
the definition of $b_c^{(n)}$ and \eqref{eq:f4o}).
Combining this with Proposition \ref{ldp:supcrit}
\begin{align}
\lim_{n\to\infty}\frac{1}{-f_4(n)\log(b_c^{(n)}/f_4(n))}\log P(B_1^{(n)})&=
\lim_{n\to\infty}\frac{a_c^{(n)}}{-f_4(n)\log(b_c^{(n)}/f_4(n))}\lim_{n\to\infty}\frac{1}{a_c^{(n)}}\log P(B_1^{(n)})\nonumber\\
&=-J(x_0)\lim_{n\to\infty}\frac{a_c^{(n)}}{-f_4(n)\log(b_c^{(n)}/f_4(n))}=-\infty.\nonumber
\end{align}
\noindent{\it Step\,\,2:\,\,Proof\,\,of\,\,\eqref{eq:compasb2terbis}.}\\
We still have the inequalities \eqref{eq:292} and so the claim follows if we check that, for $h(n)=o(p_n^{-1})$ and $\varepsilon\in\mathbb R_+$,
\begin{equation}\label{eq:limfinBISnuovo}
\lim_{n\to\infty}\frac{1}{-f_4(n)\log(b_c^{(n)}/f_4(n))}\log P(n-a_n-S_n(n-h(n))\geq\lceil\varepsilon f_4(n)\rceil)=-\varepsilon.
\end{equation}
Since $b_c^{(n)}\to 0$, by \eqref{eq:logblogbc} we have $b_c^{(n)'}\to 0$ and so by
\eqref{eq:nunomenopi} we deduce $(n-a_n)(1-\pi_n(n-h(n)))\to 0$. Consequently by \eqref{eq:nuovorate}
\begin{align}
\log P(n-a_n-S_n(n-h(n))\geq\lceil\varepsilon f_4(n)\rceil)&=\log P(\mathrm{Bin}(n-a_n,1-\pi_n(n-h(n)))\geq\lceil\varepsilon f_4(n)\rceil)\nonumber\\
&
\sim_e\varepsilon f_4(n)\log(b_c^{(n)'}/f_4(n)).\label{eq:rellogBIS}
\end{align}
By using \eqref{eq:tc}, \eqref{eq:f4o} and that $\log b_c^{(n)}\sim -n p_n$, one has $\log f_4(n)/\log b_c^{(n)}\to 0$. So, since $\log b_c^{(n)'}\sim_e\log b_c^{(n)}$,
we have
\begin{equation*}
\log(b_c^{(n)'}/f_4(n))\sim_e\log(b_c^{(n)}/f_4(n)).
\end{equation*}
Relation \eqref{eq:limfinBISnuovo} follows by \eqref{eq:rellogBIS} and this latter asymptotic equivalence.

\subsubsection{Proof of Theorem \ref{ldp3}}
\noindent{\it Proof\,\,of\,\,$(i)$.}\\
The proof is similar to the proof of Theorem \ref{ldp5} (but using Proposition \ref{ldp:supcritter}$(i)$
in place of Proposition \ref{ldp:supcritpenta}), and therefore we omit the details.\\
\noindent{\it Proof\,\,of\,\,$(ii)$.}\\
The proof is similar to the proof of Theorem \ref{ldp5} (but using Proposition \ref{ldp:supcritter}$(ii)$
in place of Proposition \ref{ldp:supcritpenta}), and therefore we omit the details.

\subsection{Proofs of Proposition \ref{ldp:supcritquater} and Theorem \ref{ldp4}}

\subsubsection{Proof of Proposition \ref{ldp:supcritquater}}

\noindent${\it Proof\,\,of\,\,(i).}$\\
Let $f_1$ be the function defined in the statement of Theorem \ref{ldp1} and let $\varepsilon\in\R_+$ be
arbitrarily fixed. By \eqref{eq:acpnk}
we have $f_1(n)=\frac{g(n)a_c^{(n)}}{n p_n}$ for $n$ large enough, and so $f_5(n)\ll f_1(n)$.
We choose $g$ in such a way that $f_1(n)=o(p_n^{-1})$ (and so, in particular, $f_1(n)/n\to 0$).
For $n\in\mathbb N$ sufficiently large, let $B_1^{(n)}$ be the event defined by \eqref{eq:b1bis} and let $B_2^{(n)}$ be the event defined by
\eqref{eq:b2ter} with $f_5$ in place of $f_4$.
We have
\begin{equation}\label{eq:maxineqquater}
P(B_1^{(n)})\vee P(B_2^{(n)})\leq P\left(\frac{n-A_n^*}{f_5(n)}>\varepsilon\right)\leq P(B_1^{(n)})+P(B_2^{(n)}).
\end{equation}
By
Proposition \ref{ldp:supcrit}
\begin{equation}\label{eq:max1}
\lim_{n\to\infty}\frac{1}{a_c^{(n)}}\log P(B_1^{(n)})=-J(x_0).
\end{equation}
Due to \eqref{eq:292} with $f_5$ in place of $f_4$, the claim then follows by the principle of the largest term (see e.g. Lemma 1.2.15 p. 7 in \cite{DZ})
if we check that, for $h(n)=o(p_n^{-1})$ and $\varepsilon\in\mathbb R_+$,
\begin{equation}\label{eq:limfin}
\lim_{n\to\infty}\frac{1}{a_c^{(n)}}\log P(n-a_n-S_n(n-h(n))\geq\lceil f_5(n)\varepsilon\rceil)=-\gamma^{-1}\lceil\ell_5\varepsilon\rceil.
\end{equation}
Note that $b_c^{(n)}\to 0$ due to \eqref{eq:acpnk} (with $\gamma\in\mathbb R_+$) and the fourth relation in \eqref{eq:ptcto0}.
By the second relation in \eqref{eq:logblogbc}
and again \eqref{eq:acpnk} (with $\gamma\in\mathbb R_+$) one has
\begin{equation}\label{eq:logasint}
\log b_c^{(n)}\sim_e -\gamma^{-1}a_c^{(n)}.
\end{equation}
Note that \eqref{eq:rellog} holds with $f_5$ in place of $f_4$. So, combining this relation with \eqref{eq:logblogbc} and \eqref{eq:logasint}, we finally have \eqref{eq:limfin}.
\\
\noindent${\it Proof\,\,of\,\,(ii).}$\\
We give the proof in the case when $b_c^{(n)}\to b:=0$ (when $b\in (0,\infty]$ the proof follows along similar arguments).
Let $f_1$, $B_1^{(n)}$ and $B_2^{(n)}$ be defined as in the proof of part $(i)$ above.
Clearly, relations \eqref{eq:maxineqquater}
and \eqref{eq:max1} still hold. Therefore,
due to \eqref{eq:292} with $f_5$ in place of $f_4$, the claim then follows by the principle of the largest term
if we check that, for $h(n)=o(p_n^{-1})$ and $\varepsilon\in\mathbb R_+$,
\begin{equation}\label{eq:limfinBIS}
\lim_{n\to\infty}\frac{1}{a_c^{(n)}}\log P(n-a_n-S_n(n-h(n))\geq\lceil\varepsilon f_5(n)\rceil)=-\ell_5'\varepsilon.
\end{equation}
To this aim, we start noticing that it is easily realized that
\eqref{eq:rellogBIS} holds with $f_5$ in place of $f_4$. Therefore, by the assumption
$f_5(n)\sim_e\ell_5'\frac{a_c^{(n)}}{n p_n}$,
\[
(a_c^{(n)})^{-1}\log P(n-a_n-S_n(n-h(n))\geq\lceil\varepsilon f_5(n)\rceil)
\sim_e\ell_5'\varepsilon\frac{\log((n p_n b_c^{(n)'})/a_c^{(n)})}{n p_n}.
\]
Relation \eqref{eq:limfinBIS} follows noticing that
\begin{align}
\frac{\log((n p_n b_c^{(n)'})/a_c^{(n)})}{n p_n}&=\frac{\log(n p_n)}{n p_n}+\frac{\log b_c^{(n)'}}{n p_n}-\frac{\log a_c^{(n)}}{n p_n}\nonumber\\
&\sim_e -1-\frac{\log a_c^{(n)}}{n p_n}\sim_e -1,\label{eq:maggio1}
\end{align}
where the first asymptotic equivalence in \eqref{eq:maggio1} is a consequence of $n p_n\to\infty$ and \eqref{eq:logblogbc}, and
the second asymptotic equivalence in \eqref{eq:maggio1} follows by
$(\log a_c^{(n)})/(n p_n)\to 0$ (which is an easy consequence of the definition of $a_c^{(n)}$).

\subsubsection{Proof of Theorem \ref{ldp4}}

The proof is very similar to the proof of Theorem \ref{ldp1}
and therefore we omit the details (clearly one has to
exploit Proposition \ref{ldp:supcritquater} in place of Proposition \ref{ldp:supcrit}).

\section{Appendix}\label{sec:app}

\subsection{Proofs of \eqref{eq:Jasonbin} and \eqref{eq:nuovorate}}

The first part of the proof is common for \eqref{eq:Jasonbin} and \eqref{eq:nuovorate} (clearly, dealing with \eqref{eq:Jasonbin} we have to set $r_n=k\in\mathbb N$ for any $n$).
By assumption $m_n/r_n\to\infty$ and so, for all $n$ large enough,
\begin{align}
P(\mathrm{Bin}(\lfloor m_n\rfloor,q_n)\geq r_n)&=\sum_{k=r_n}^{\lfloor m_n\rfloor}\binom{\lfloor m_n\rfloor}{k}q_n^{k}(1-q_n)^{\lfloor m_n\rfloor-k}\nonumber\\
&=\binom{\lfloor m_n\rfloor}{r_n}q_n^{r_n}(1-q_n)^{\lfloor m_n\rfloor-r_n}
\sum_{k=r_n}^{\lfloor m_n\rfloor}\binom{\lfloor m_n\rfloor}{r_n}^{-1}
\binom{\lfloor m_n\rfloor}{k}q_n^{k-r_n}(1-q_n)^{r_n-k}\nonumber\\
&=\binom{\lfloor m_n\rfloor}{r_n}q_n^{r_n}(1-q_n)^{\lfloor m_n\rfloor-r_n}
\sum_{j=0}^{\lfloor m_n\rfloor-r_n}\binom{\lfloor m_n\rfloor}{r_n}^{-1}
\binom{\lfloor m_n\rfloor}{r_n+j}\left(\frac{q_n}{1-q_n}\right)^{j}\nonumber\\
&\leq\binom{\lfloor m_n\rfloor}{r_n}q_n^{r_n}(1-q_n)^{\lfloor m_n\rfloor-r_n}
\sum_{j=0}^{\lfloor m_n\rfloor-r_n}\left(\frac{\lfloor m_n\rfloor q_n}{r_n(1-q_n)}\right)^{j}.\nonumber
\end{align}
By assumption $q_n\to 0$ and $(m_n q_n)/r_n\to 0$, and so (as it may be easily checked)
\[
\sum_{j=0}^{\lfloor m_n\rfloor-r_n}\left(\frac{\lfloor m_n\rfloor q_n}{r_n(1-q_n)}\right)^{j}\to 1.
\]
Consequently
\[
P(\mathrm{Bin}(\lfloor m_n\rfloor,q_n)\geq r_n)\sim_e\binom{\lfloor m_n\rfloor}{r_n}q_n^{r_n}(1-q_n)^{\lfloor m_n\rfloor-r_n}.
\]
Since $m_n\to\infty$ and $m_n-r_n\to\infty$, by Stirling's formula
\begin{align}
\frac{\lfloor m_n\rfloor!}{(\lfloor m_n\rfloor-r_n)!}&\sim_e\sqrt{\frac{m_n}{m_n-r_n}}\left(\frac{m_n}{\mathrm{e}}\right)^{r_n}\left(\frac{m_n}{m_n-r_n}\right)^{m_n-r_n}\nonumber\\
&\sim_e\sqrt{\frac{m_n}{m_n-r_n}}\left(\frac{m_n}{\mathrm{e}}\right)^{r_n}\mathrm{e}^{r_n}\sim_e m_n^{r_n},\nonumber
\end{align}
and so
\begin{equation}\label{eq:binasrel}
P(\mathrm{Bin}(\lfloor m_n\rfloor,q_n)\geq r_n)\sim_e\frac{(m_n q_n)^{r_n}}{r_n!}(1-q_n)^{m_n-r_n}.
\end{equation}
We proceed by distinguishing the proof of \eqref{eq:Jasonbin} and the proof of \eqref{eq:nuovorate}.

\subsubsection{Proof of \eqref{eq:Jasonbin}}

If $r_n=k\in\mathbb N$ for any $n$, then \eqref{eq:Jasonbin} easily follows by \eqref{eq:binasrel}. Indeed
\[
(1-q_n)^{m_n-k}=\left((1-q_n)^{q_n^{-1}}\right)^{q_n(m_n-k)}\to\mathrm{e}^{0}=1,
\]
where the limit follows by $q_n\to 0$ and $q_n m_n\to 0$.

\subsubsection{Proof of \eqref{eq:nuovorate}}

If $r_n\to\infty$, then by Stirling's formula
\[
r_n!\sim_e\sqrt{2\pi r_n}(r_n/\mathrm{e})^{r_n},
\]
and so \eqref{eq:binasrel} yields
\begin{equation*}
P(\mathrm{Bin}(\lfloor m_n\rfloor,q_n)\geq r_n)\sim_e\frac{(m_n q_n\mathrm{e}/r_n)^{r_n}}{\sqrt{2\pi r_n}}(1-q_n)^{m_n-r_n}.
\end{equation*}
Consequently,
\begin{align}
\log P(\mathrm{Bin}(\lfloor m_n\rfloor,q_n)\geq r_n)&\sim_e r_n\log\left(\frac{m_n q_n\mathrm{e}}{r_n}\right)-\frac{\log r_n}{2}+(m_n-r_n)\log(1-q_n)\nonumber\\
&=r_n+r_n\log\left(\frac{m_n q_n}{r_n}\right)-\frac{\log r_n}{2}+(m_n-r_n)\log(1-q_n)\nonumber\\
&=r_n\log\left(\frac{m_n q_n}{r_n}\right)\Biggl[\frac{1}{\log\left(\frac{m_n q_n}{r_n}\right)}-\frac{\log r_n}{2 r_n\log\left(\frac{m_n q_n}{r_n}\right)}+1\nonumber\\
&\,\,\,\,\,\,
\,\,\,\,\,\,\,
\,\,\,\,\,\,\,
\,\,\,\,\,\,\,
\,\,\,\,\,\,\,
\,\,\,\,\,\,\,
\,\,\,\,\,\,\,
\,\,\,\,\,\,\,
\,\,\,\,\,\,\,
\,\,\,\,\,\,\,
\,\,\,\,\,\,\,
+\frac{(m_n-r_n)\log(1-q_n)}{r_n\log\left(\frac{m_n q_n}{r_n}\right)}\Biggr].\nonumber
\end{align}
Since $m_n q_n/r_n\to 0$ the first two addends of the quantity between the squared brackets tend to $0$. Exploiting also that $m_n/r_n\to\infty$ and $q_n\to 0$, we have
\begin{align}
\frac{(m_n-r_n)\log(1-q_n)}{r_n\log\left(\frac{m_n q_n}{r_n}\right)}&\sim_e -\frac{(m_n q_n/r_n)}{\log(m_n q_n/r_n)}\to 0,\nonumber
\end{align}
which concludes the proof.

\subsection{Proofs of \eqref{eq:logblogbc1} and \eqref{eq:logblogbc}}

\subsubsection{Proof of \eqref{eq:logblogbc1}}

If $b_c^{(n)}\to b:=\infty$, then combining formulas $(3.7)$ and $(3.9)$ in \cite{JLTV}
one has $b_c^{(n)}\sim_e b_c^{(n)'}$. If $b_c^{(n)}\to b\in\mathbb R_+$, then by \eqref{rfIexpression}
$n p_n-(\log n +(r-1)\log\log n)\to-\log((r-1)!b)$, and so $p_n\sim_e\log n/n$. Consequently, $p_n=o(n^{-1/2})$ and by formula $(3.8)$
in \cite{JLTV} again $b_c^{(n)}\sim_e b_c^{(n)'}$.

\subsubsection{Proof of \eqref{eq:logblogbc}}

By \eqref{rfIexpression} one has that there exists $M>1$ and $n_0\in\mathbb N$: $n p_n\geq M\log n$ for all $n\geq n_0$.
Therefore $\log n/(n p_n)\to\mathrm{const}$, where either $\mathrm{const}=M'\in (0,M^{-1}]$ or $\mathrm{const}=0$. Consequently
\begin{align}
\lim_{n\to\infty}\frac{\log b_c^{(n)}}{\log b_c^{(n)'}}&=\lim_{n\to\infty}\frac{\log n+(r-1)\log(n p_n)-n p_n}{\log n+(r-1)\log(n p_n)+n\log(1-p_n)}\nonumber\\
&=\lim_{n\to\infty}\frac{\log n/(n p_n)+(r-1)\log(n p_n)/(n p_n)-1}{\log n/(n p_n)+(r-1)\log(n p_n)/(n p_n)+\log(1-p_n)/p_n}=1\nonumber
\end{align}
and
\begin{align}
\lim_{n\to\infty}\frac{\log b_c^{(n)}}{-np_n}&=\lim_{n\to\infty}\frac{\log n+(r-1)\log(n p_n)-n p_n}{-n p_n}=-M'+1>0.\nonumber
\end{align}

\end{document}